\documentclass[10pt,twoside]{article}
\usepackage{amsmath,amsthm}
\usepackage[latin1]{inputenc}
\usepackage{picinpar}
\usepackage{longtable}
\usepackage{amsrefs}
\usepackage[dvips]{graphicx}
\usepackage[mathscr]{eucal}
\usepackage{calc}
\usepackage{ifthen}
\usepackage{dcpic,pictexwd}
\usepackage[all]{xy}
\usepackage{enumerate}
\usepackage{amssymb,latexsym}
\usepackage{amsthm}
\usepackage[dvips]{graphics}
\usepackage{color}
\usepackage{tabularx}
\usepackage{makeidx}
\usepackage{mathdots}
\usepackage{tikz-cd}
\usepackage[ruled,vlined]{algorithm2e}
\usetikzlibrary{calc}
\usepackage{tikz}
\usetikzlibrary{arrows}
\usepackage{stackrel}
\usepackage{multicol}
\usepackage{float}
\newtheorem{teorema}{Theorem}

\newtheorem{corolario}{Corollary}

\newtheorem{proposicion}{Proposition}

 \newtheorem{teor}{Theorem}

\newtheorem{prop}{Proposition}
\newtheorem{corol}{Corollary}

\newtheorem{Nota}{Remark}

\usepackage{textcomp}
\usepackage{tikz}
\theoremstyle{definition}
\theoremstyle{definition}

\renewcommand{\figurename}{Figure}

\begin{document}
\renewcommand{\figurename}{Fig.}

\thispagestyle{plain}
\par\bigskip
\begin{centering}

\textbf{ Relationships Between Mutations of Brauer Configuration Algebras and Some Diophantine Equations}

\end{centering}\par\bigskip

\par\bigskip
\begin{centering}

\footnotesize{Agustín Moreno Ca\~{n}adas}\\
\footnotesize{Juan David Camacho}\\
\footnotesize{Isa\'{i}as David Mar\'{i}n Gaviria}

\end{centering}\par\bigskip

\bigskip\begin{centering} \textbf{Abstract}\par\bigskip\end{centering} Mutations on Brauer configurations are introduced and associated with some suitable automata in order to solve generalizations of the Chicken McNugget problem. Besides, based on marked order polytopes a new class of diophantine equations called Gelfand-Tsetlin equations are also solved. The approach allows giving an algebraic description of the schedule of an AES key via some suitable non-deterministic finite automata (NFA).


\bigskip

\section{Introduction}
Brauer configuration algebras were introduced recently by Green and Schroll \cite{Green}. They constitute an amazing mathematical object with applications in a wide variety of fields. For instance, Ca\~{n}adas et al. have used Brauer configuration algebras to describe
secret sharing schemes in order to protect biometric data \cite{MAO}. On the other hand, Fern\'andez et al defined the message associated with a Brauer configuration in order to interpret different mathematical notions as a suitable message. Snake graphs, homological ideals associated with algebras of Dynkin type $\mathbb{A}_{n}$ and Dyck paths are examples of different mathematical objects whose enumeration can be obtained via messages of Brauer configuration algebras \cites{Fernandez, Fernandez1}.\par\bigskip

Brauer configuration algebras arose as a generalization of Brauer graph algebras whose theory of representation is given by  combinatorial data. Perhaps this is one of the reasons for which these algebras appear in many fields in the form of messages. In this paper, we introduce the notion of mutation associated with Brauer configurations and the messages of these mutations in order to solve generalizations of the so-called Chicken McNugget Problem, which arose from the way that the international fast food restaurant chain McDonald\textquotesingle s delivers one of its most sold product, the Chicken McNuggets.\par\bigskip

We recall that shortly after their introduction in 1979, Chicken McNuggets were sold in packs of 6, 9, and 20 pieces. Nowadays, several choices of this product can be ordered, which in a natural fashion generate different types of diophantine equations. For instance, the first three types of packs give rise to the following problem \cite{Chapman}:\par\bigskip
\begin{centering}
\textbf{Problem A.}
\end{centering}
\textit{What numbers of chicken McNuggets can be ordered using only packs with 6, 9, or 20 pieces}?\par\bigskip

Problem A is currently called the \textit{Chicken McNugget Problem} (CMP), which defines the following diophantine equation:

\begin{equation}\label{McNuggets equation}
\begin{split}
6x+9y+20z&=n
\end{split}
\end{equation}
 where $n$ is the number of ordered pieces.
\par\bigskip

CMP is a version of the \textit{Postage Stamp Problem} or \textit{Frobenius Coin Problem} or \textit{Knapsack Problem},  which can be defined as follows:\par\bigskip

\textbf{Problem B.} \textit{Given a set of $k$ objects with predetermined values $n_{1}, n_{2},\dots, n_{k}$, what possible values of $n$ can be had from combinations of these objects}?\par\bigskip
 
Positive integers satisfying the CMP are known as McNugget numbers, sequence A065003 in the OEIS is the list of not McNugget numbers. Actually, it is known that any positive integer $x>43$ is a McNugget number \cite{Chapman}.\par\bigskip

CMP deals with a more wide kind of problems based on the number of lattice points in a convex polytope $P$, i.e., to determine $\phi_{A}(b)=|P=\{x\mid Ax=b,\hspace{0.1cm}x\hspace{0.1cm}\text{integral}\}|$, where $A=(a_{(i,j)})$ is a $m\times n$ integral matrix and $b$ is an integral $m$-vector. For instance, if $H=[3\hspace{0.2cm}5\hspace{0.2cm}17]$ then
$\phi_{H}(58)=9$, $\phi_{H}(101)=25$, $\phi_{H}(1110)=2471$ \cite{DeLoera}. Actually,
\begin{equation}
 \underset{n=0}{\overset{\infty}{\sum}}\phi_{H}(n)t^{n}=\frac{1}{(1-t^3)(1-t^{5})(1-t^{17})}\notag
\end{equation}

in general it holds that:

\begin{equation}
 \underset{n=0}{\overset{\infty}{\sum}}\phi_{A}(n)t^{n}=\frac{1}{1-t^{a_{(1,j)}a_{(2,j)}a_{(3,j)}\dots a_{(m, j)}}}.\notag
\end{equation}
\par\bigskip
It is worth noting that in general solving an equation as the aforementioned (of type $Ax=b$) is a cumbersome problem. 
In \cite{DeLoera} De Loera gives examples of hard instances of this type. \par\bigskip

 Gelfand-Tsetlin polytopes (or GT-polytopes) are another kind of polytopes whose enumeration gives rise to interesting problems in combinatorics, diophantine analysis and many other fields. 
We recall that a $T$-array, $T=(t_{(i,j)})\in\mathbb{Z}^{\frac{n(n+1)}{2}}$ is called a \textit{standard Gelfand-Tsetlin  pattern} (GT-pattern or Gelfand-Tsetlin tableau), if it satisfies the interlacing conditions \cite{Futorny}:
\begin{equation}\label{GT1}
\begin{split}
t_{(i+1,j)}-t_{(i,j)}&\geq0,\\
t_{(i,j)}-t_{(i+1,j+1)}>0.
\end{split}
\end{equation}

Given $\lambda,\mu\in\mathbb{Z}^{n}$, the \textit{Gelfand-Tsetlin polytope} $\mathrm{GT}(\lambda,\mu)$ is the convex polytope of $\mathrm{GT}$-patterns $(t_{(i, j)})_{1\leq i\leq j\leq n}$ satisfying in addition that

\begin{equation}\label{GT2}
\begin{split}
t_{(1,1)}&=\mu_{1},\\
\underset{i=1}{\overset{j}{\sum}}t_{(i,j)}-\underset{i=1}{\overset{j-1}{\sum}}t_{(i,j-1)}&=\mu_{j},\quad 2\leq j\leq n.
\end{split}
\end{equation}

It should be pointed out that much of the theory of representation of the Lie algebra $\mathrm{gl}_{n}\hspace{0.1cm}\mathbb{C}$ is based on a result of Gelfand and Tsetlin related with the enumeration of lattice points in GT-polytopes \cites{DeLoera, Futorny, RamirezF}. Actually, it is known that if $U$ is the subspace of an irreducible representation $V$ of $\mathrm{gl}_{n}\hspace{0.1cm}\mathbb{C}$ with highest weight $\lambda$ and $d=\mathrm{dim}_{\mathbb{C}}\hspace{0.1cm}U$ then $d$ is given by the number of integral points in $\mathrm{GT}(\lambda,\mu)$.\par\bigskip

In this work, we introduce Gelfand-Tsetlin numbers that satisfy some diophantine equations called Gelfand-Tsetlin equations. In particular, we define Brauer configurations whose vertices are lattice points (of a lattice polytope $\mathscr{P}_{(n, r)}$) in bijective correspondence with sets of Gelfand-Tsetlin patterns. Brauer configurations also allow us to define NFAs whose states are given by solutions of diophantine problems of the form $\mathcal{D}(n_{1},n_{2}, \mathcal{K}_{m})$, in this case, $n_{1}\leq 16$, $\mathcal{K}_{m}=\{k_{1},k_{2},\dots, k_{m}\}$ is a fixed set of occurrences of vertices of some suitable Brauer configurations, and for the fixed integers $n_{1}$ and $n_{2}$, it holds that:

\begin{equation}\label{type D}
\begin{split}
\underset{i=1}{\overset{m}{\sum}}\lambda_{i}&=n_{1},\\
\underset{j=1}{\overset{m}{\sum}}k_{j}\lambda_{j}&=n_{2}.
\end{split}
\end{equation}

This approach permits to define the cryptographic AES-key schedule as the message of a Brauer configuration algebra. To do that, the notion of mutation of a Brauer configuration is introduced.
\par\bigskip

\subsection{Main Results}

In the sequel, we give a brief description of some of the main results presented in this paper.\par\bigskip

The following result gives some properties of some Brauer configuration algebras $\Lambda_{\Phi^{m_{0}}}$ obtained upon a mutation process is completed (see Theorem \ref{MUT}).
\addtocounter{teorema}{2}
\begin{teorema}\label{MutA}
For $m_{0}>1$ the Brauer configuration algebra $\Lambda_{\Phi^{m_{0}}}$ induced by the Brauer configuration $\Phi^{m_{0}}$ (see $\mathrm{(\ref{mutation})}$) is connected and reduced.
\end{teorema}

Regarding, enumeration of some Gelfand-Tsetlin arrays, we give the following results:

\begin{teorema}\label{NGTB}
The number $g_{(n, r)}$ of GT-patterns over $\mathrm{gl}_{n}$ defined by a weight vector of the form
$w=(n, n-r, n-2r,\dots, n-r(n-1))$ is given by the formula $g_{(n, r)}=(r+1)^{\binom{n}{2}}$.
\end{teorema}

In the next result, $\mathscr{H}_{n}(\tau)$ denotes a set of Gelfand-Tsetlin arrays $\tau$ associated with a suitable triangular array called Gelfand-Tsetlin heart, whereas $T(\mathcal{H}_{n})$ denotes a set of suitable Gelfand-Tsetlin arrays with $\mathcal{H}_{n}$ as its heart. 
\addtocounter{proposicion}{4}
\begin{proposicion}\label{NGTC}
\textit{If $n=4$, $r>1$, and $\mathcal{H}_{n},\mathcal{H}'_{n}\in\mathscr{H}'_{(n, r)}$ are such that $\lambda_{(n-1,2)}=\lambda'_{(n-1,2)}$ and
$\lambda_{(n-2,1)}=\lambda'_{(n-2,1)}$ then $|T(\mathcal{H}_{n})|=|T(\mathcal{H}'_{n})|$.}

\end{proposicion}

In the following result $Q_{\mathscr{H}_{(n, r)}}$ is a notation for some suitable marked posets.

\begin{proposicion}\label{NGT2D}
\textit{For $n=4$, and $r>1$, if $\mathcal{Q}$ denotes the set of numbers $q_{i}$ given by identities $\mathrm{(\ref{sequence})}$ and $\mathrm{(\ref{sequence1})}$. Then $\mathcal{Q}\subset Q_{\mathscr{H}_{(n, r)}}$ (see $\mathrm{(\ref{markedposet})}$).}

\end{proposicion}

\addtocounter{corolario}{6}

\begin{corolario}
For $r\geq 2$ and $n=4$, the number of facets $\mathscr{F}_{(n, r)}$ in the marked order polytope $Q_{\mathscr{H}_{(n, r)}}$ associated with $\mathscr{H}_{(n, r)}$ is $r(r+1)(3r+2).$
\end{corolario}

The following result regards Gelfand-Tsetlin numbers $S_{gt}$, which are solutions of Gelfand-Tsetlin equations.

\begin{corolario}
\textit{For $r>1$ and $n=4$ sums  $S_{gt}(n_{1},r)$ (see identities $\mathrm{(\ref{sequence})}$ and $\mathrm{(\ref{sequence1})}$) associated with points $\mathcal{H}_{n}\in(\mathscr{H}_{(n,r)},\unlhd)$ are Gelfand-Tsetlin numbers of type $(n,r)$.}
\end{corolario}

The following results give relationships between Brauer configuration algebras and diophantine equations.

\addtocounter{teorema}{4}

\begin{teorema}
If $\mathbb{F}$ is a finite field then the set $M_{\Phi}$ of $m_{0}$-Brauer clusters obtained by mutation is finite, and the special cycle of maximal length has one of the two shapes shown below. Moreover, if $\Phi^{i}\neq\Phi^{j}$ whenever $i\neq j$ then there exists $m\in\mathbb{N}$ such that  $(\Gamma, \mathscr{X})= (\Gamma^{m+1}, \mathscr{X}^{m+1})$ for a given initial seed  $(\Gamma, \mathscr{X})$.

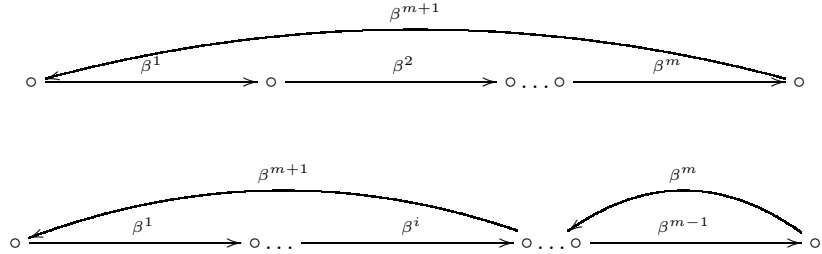
\begin{equation}\label{mut1s}
\xymatrix@=80pt{{\circ}\ar@/^0pt/@[\blue][r]^{\beta^{1}}&{\circ}\ar@/^0pt/@[\red][r]^{\beta^2}& \circ\dots \ar@/^0pt/@[\red][r]^{\beta^m}\circ&\circ \ar@/_20pt/@[\red][lll]_{\beta^{m+1}}}\notag
\end{equation}

\begin{equation}\label{mut1sa}
\xymatrix@=80pt{{\circ}\ar@/^0pt/@[\blue][r]^{\beta^{1}}&{\circ}\dots\ar@/^0pt/@[\red][r]^{\beta^i}& \circ\dots \ar@/^0pt/@[\red][r]^{\beta^{m-1}}\ar@/_20pt/@[\red][ll]_{\beta^{m+1}}\circ&\circ \ar@/_20pt/@[\red][l]_{\beta^{m}}}\notag
\end{equation}

where $\beta^{m}\beta^{m+1}\neq0$.

\end{teorema}

\begin{teorema}\label{mainH}
For a fixed integer positive $m_{0}$, $n\geq2$, and a fixed seed $(\Gamma,\mathscr{X})$ (see $\mathrm{(\ref{mut01})}$ and $\mathrm{(\ref{mut02})}$), the Brauer configuration algebra $\Lambda_{\Phi^{m_{0}}}$ induced by the Brauer configuration 
\begin{equation}
\Phi^{m_{0}}=(\Phi^{m_{0}}_{0},\Phi^{m_{0}}_{1},\mu,\mathcal{O})\quad (see\hspace{0.1cm} \mathrm(\ref{mutation}))\notag
\end{equation}

has associated a finite non-deterministic automaton whose states are given by solutions of a problem of type $\mathcal{D}(n_{1},l^{2}2^{n-2},\mathcal{K}_{m_{i}})$ with $n_{1}\leq 16$, and 
\begin{equation}
\mathcal{K}_{m_{i}}= \{\mathrm{occ}(T(\alpha), T(\mathcal{M}^{i}))\mid T(\alpha)\in T(\Phi^{m_{0}}_{0}), T(\mathcal{M}^{i})\in T(\Phi^{m_{0}}_{1}), \alpha\in\Phi^{m_{0}}_{0}, \mathcal{M}^{i}\in \Phi^{m_{0}}_{1} \}\notag
\end{equation}

where $T(\Phi^{m_{0}})=((T(\Phi^{m_{0}}_{0}),T(\Phi^{m_{0}}_{1}),T(\mu),T(\mathcal{O}))$ is a suitable transformation between Brauer configurations.
\end{teorema}

\par\bigskip
This paper is organized as follows; Basic notation and definitions to be used throughout the paper are given in section 2. In particular, we introduce an algorithm to build Brauer configuration algebras, and the notion of mutation associated with Brauer configurations is also introduced in section 2.2.2. In section 3, we prove results regarding Gelfand-Tsetlin patterns, in particular, notions of Gelfand-Tsetlin equations and Gelfand-Tsetlin numbers are introduced in this section as well. In section 4, we give properties of some diophantine equations whose solutions are given by mutations of some Brauer configurations. In this section, the schedule of an AES key is interpreted as a mutation process associated with suitable Brauer configurations.

\section{Preliminaries}

In this section, we introduce basic definitions and notation
to be used throughout the paper. Henceforth, we will use the customary symbols $\mathbb{N}$, $\mathbb{Z}$,
$\mathbb{R}$, and $\mathbb{C}$ to denote the set of natural numbers, integers, real, and complex
numbers respectively.

 \subsection{On the Frobenius Number}

 If $n_{1}, n_{2},\dots, n_{k}$ are positive integers then the set $S$ of all integers which can be presented in the form \cite{Chapman}:
 
 \begin{equation}
 \begin{split}
n&=\underset{i=1}{\overset{k}{\sum}}n_{i}\lambda_{i},\quad \lambda_{1},\lambda_{2},\dots,\lambda_{k}\in\mathbb{N} 
 \end{split}
 \end{equation}
 
is a submonoid of $(\mathbb{N},+)$. $S$ is said to be a \textit{numerical semigroup} if $|\mathbb{N}\backslash S|<\infty$.\par\bigskip 

$\langle n_{1}, n_{2},\dots, n_{k}\rangle=\{\underset{i=1}{\overset{k}{\sum}}\lambda_{i}n_{i}\mid \lambda_{i}\in\mathbb{N}, 1\leq i\leq k\}$ is said to be a \textit{numerical monoid} generated by $\{n_{1},n_{2},\dots, n_{k}\}$.  $\langle 6, 9, 20\rangle$ is the \textit{chicken McNugget monoid}.\par\bigskip

The following result regards numerical semigroups.

\begin{teor}[\cite{Garcia}, p.7]
\textit{If $n_{1},n_{2},\dots, n_{k}$ are positive integers, then they generate a numerical semigroup if and only if $(n_{1},n_{2},\dots, n_{k})=1$.}
\end{teor}

 If $n_{1},n_{2},\dots, n_{k}$ are positive integers, then the \textit{Frobenius number} of $\langle n_{1},n_{2},\dots, n_{k}\rangle$ denoted $F(\langle n_{1},n_{2},\dots, n_{k}\rangle)$ is the largest positive integer $n$ such that $n\notin \langle n_{1},n_{2},\dots, n_{k}\rangle$. For instance, $F(\langle 6, 9, 20,\rangle)=43$. In the general case, is not known a formula for the Frobenius number of the $k$-generated numerical semigroup $S$, but for $k$ fixed, there is an algorithm that computes the Frobenius number in polynomial time. Actually, Ram\'{\i}rez-Alfons\'{\i}n proved that the knapsack problem can be reduced to the Frobenius problem in polynomial time \cite{Ramirez}. Let $x, y\in \langle n_{1},n_{2},\dots, n_{k}\rangle$, we say that \textit{$x$ divides $y$} in $\langle n_{1},n_{2},\dots, n_{k}\rangle$ if there exists $z\in\langle n_{1},n_{2},\dots, n_{k}\rangle$ such that $y=x+z$. 
 
\par\bigskip 
 We call a nonzero element $x\in\langle n_{1},n_{2},\dots, n_{k}\rangle$ \textit{irreducible} if whenever $x=y+z$, either $y=0$ or $z=0$ (hence, $x$ is irreducible if its only proper divisors are only 0 or itself).
\par\bigskip
The following result determines which elements are irreducible in a numerical monoid.

\addtocounter{prop}{1}
\begin{prop}[\cite{Chapman}, Proposition 3.1]
\textit{If $\langle n_{1},n_{2},\dots, n_{k}\rangle$ is a numerical monoid, then its irreducible elements are precisely $n_{1}, n_{2},\dots, n_{k}$. }

\end{prop}

 \addtocounter{corol}{2}
 \begin{corol}[\cite{Chapman}, Corollary 3.2]
 \textit{The irreducible elements of the McNugget monoid are 6, 9, and 20.}
 \end{corol}

 \subsection{Brauer Configuration Algebras}

As we mentioned in the introduction, Brauer configuration algebras are multiserial algebras introduced recently by Green and Schroll in \cite{Green}. These algebras constitute a generalization of Brauer graph algebras, which have as one of their properties that its theory of representation is encoded by some combinatorial data based on graphs. 

\par\bigskip

The following is a description of the structure of Brauer configuration algebras as the first author and Espinosa and Green and Schroll present in \cite{Fernandez1} and \cite{Green}, respectively.\par\bigskip

A \textit{Brauer configuration} $\Gamma$ is a quadruple of the form $\Gamma=(\Gamma_{0},\Gamma_{1},\mu,\mathcal{O})$, where

\begin{enumerate}
\item $\Gamma_{0}$ is a set whose elements are called \textit{vertices}, $|\Gamma_{0}|>1$,
\item $\Gamma_{1}$ is a finite collection of multisets called polygons, such that if $V\in\Gamma_{1}$ then $|V|>1$,
\item $\mu$ is an integer valued map, $\mu:\Gamma_{0}\rightarrow \mathbb{N}\backslash\{0\}$ (i.e., for each vertex $\alpha\in\Gamma_{0}$, it holds that $\mu(\alpha)>1$) called the \textit{multiplicity function}.
\item $\mathcal{O}$ is a well-ordering $<$ defined on $\Gamma_{1}$. Such that for each vertex $\alpha\in\Gamma_{0}$, the collection $S_{\alpha}=\{V^{(\alpha_{1})}_{\alpha,1}, V^{(\alpha_{2})}_{\alpha,2},\dots, V^{(\alpha_{t})}_{\alpha,t}\mid \alpha_{i}\geq1\hspace{0.1cm}\text{denotes the frequency of}\hspace{0.1cm}\alpha\hspace{0.1cm}\text{in}\hspace{0.1cm}V_{\alpha,i}\}$ consisting of all polygons of $\Gamma_{1}$ where $\alpha$ occurs (counting repetitions) is ordered in the form $V^{(\alpha_{1})}_{\alpha,1}<V^{(\alpha_{2})}_{\alpha,2}<\dots< V^{(\alpha_{t})}_{\alpha,t}$ ($V_{\alpha,i}<V_{\alpha,i}$, for $1\leq i\leq t$). If it is assumed that $V_{\alpha,1}=\mathrm{min}\hspace{0.1cm}S_{\alpha}$ and $V_{\alpha,t}=\mathrm{max}\hspace{0.1cm}S_{\alpha}$ then a new circular relation $R_{\alpha}=\{V_{\alpha,t}<V_{\alpha,1}\}$ is added. $S_{\alpha}$ is called the \textit{successor sequence} at the vertex $\alpha$ and $C_{\alpha}=S_{\alpha}\cup R_{\alpha}$.
\item $\mathrm{occ}(\alpha, V)$ denotes the number of times that a vertex $\alpha$ occurs in a polygon $V$. And the sum $\underset{V\in\Gamma_{1}}{\sum}\mathrm{occ}(\alpha, V)$ is said to be the \textit{valency} of $\alpha$, denoted $val(\alpha)$.

\item Vertices are either \textit{truncated} or \textit{non-truncated} (i.e., $\Gamma_{0}=T_{\Gamma}\overset{.}{\cup}NT_{\Gamma}$, where $T_{\Gamma}$ ($NT_{\Gamma}$) is the set of truncated (non-truncated) vertices), if $U\in\Gamma_{1}$ then $U\cap NT_{\Gamma}\neq\varnothing$. A vertex $\alpha\in\Gamma_{0}$ is truncated if $\mu(\alpha)val(\alpha)=1$, it is non-truncated if $\mu(\alpha)val(\alpha)>1$. A Brauer configuration $\Gamma$ without truncated vertices is \textit{reduced}, these types of configurations are used to construct suitable quivers called \textit{Brauer quivers}.

\end{enumerate}

The following algorithm builds the Brauer quiver $Q_{\Gamma}$ and the Brauer configuration algebra $\Lambda_{\Gamma}=\mathbb{F}Q_{\Gamma}/I_{\Gamma}$ induced by the Brauer configuration $\Gamma$, where $I_{\Gamma}$ is an admissible ideal (see Remark \ref{class}).
\par\bigskip

\begin{algorithm}[H]\label{algoritmo0}
\SetAlgoLined
\begin{enumerate}
\item \textbf{Input} A reduced Brauer configuration $\Gamma=(\Gamma_{0},\Gamma_{1},\mu,\mathcal{O})$.
\item\textbf{Output} The Brauer configuration algebra $\Lambda_{\Gamma}=\mathbb{F}Q_{\Gamma}/I_{\Gamma}$.
\item Construct the quiver $Q_{\Gamma}=((Q_{\Gamma})_{0},(Q_{\Gamma})_{1},s:(Q_{\Gamma})_{1}\rightarrow(Q_{\Gamma})_{0},t:(Q_{\Gamma})_{1}\rightarrow(Q_{\Gamma})_{0})$
\begin{enumerate}
\item $(Q_{\Gamma})_{0}=\Gamma_{1}$,
\item For each cover $V_{i}<V_{i+1}\in\Gamma_{1}$ define an arrow $a\in(Q_{\Gamma})_{1}$, such that $s(a)=V_{i}$ and $t(a)=V_{i+1}$,
\item Each relation $V_{i}<V_{i}$ defines a loop in $Q_{\Gamma}$,
\item Each ordered set $C_{\alpha}$ defines a cycle in $Q_{\Gamma}$ called \textit{special cycle}.
\end{enumerate}

\item Define the path algebra $\mathbb{F}Q_{\Gamma}$.

\item Construct $I_{\Gamma}$, which is generated by the following relations:
\begin{enumerate}
\item If $\alpha_{i},\alpha_{j}\in U$, $U\in\Gamma_{1}$ and $C_{\alpha_{i}}, C_{\alpha_{j}}$ are corresponding special cycles then $C^{\mu(\alpha_{i})}_{\alpha_{i}}-C^{\mu(\alpha_{j})}_{\alpha_{j}}=0$,
\item  If $C_{\alpha_{i}}$ is a special cycle associated with the vertex $\alpha_{i}$ then $C^{\mu(\alpha_{i})}a=0$, if $a$ is the first arrow of
$C_{\alpha_{i}}$,
\item If $\alpha,\alpha'\in\Gamma_{0}$, $\alpha\neq\alpha'$, $a\in C_{\alpha}$, $b\in C_{\alpha'}$, $a\neq b$, $ab\notin C_{\alpha}$ for any $\alpha\in\Gamma_{0}$ then $ab=0$,
\item If $a$ is a loop associated with a vertex $\alpha$ with $val(\alpha)=1$ and $\mu(\alpha)>1$ then $a^{\mu(\alpha)+1}=0$.
\end{enumerate}
\item $\Lambda_{\Gamma}=\mathbb{F}Q_{\Gamma}/I_{\Gamma}$ is the Brauer configuration algebra. 
\item For the construction of a basis of $\Lambda_{\Gamma}$ follow the next steps:
\begin{enumerate}
\item For each $V\in\Gamma_{1}$ choose a non-truncated vertex $\alpha_{V}$ and exactly one special $\alpha$-cycle $C_{\alpha_{V}}$ at $V$,
\item Define: 
\begin{equation*}
\begin{split}
A&=\{\overline{p}\mid p\hspace{0.1cm} is\hspace{0.1cm} a\hspace{0.1cm} proper\hspace{0.1cm} prefix\hspace{0.1cm} of\hspace{0.1cm} some\hspace{0.1cm}C^{\mu(\alpha)}_{\alpha}\},\\
B&= \{\overline{C^{\mu(\alpha)}_{\alpha_{V}}}\mid V\in \Gamma_{1}\}.
\end{split}
\end{equation*}
\item $A\cup B$ is an $\mathbb{F}$-basis of $\Lambda_{\Gamma}$.
\end{enumerate}

\end{enumerate}
\caption{Construction of a Brauer configuration algebra}
\end{algorithm}

\addtocounter{Nota}{3}
\begin{Nota}\label{class}
Let $\Lambda=\mathbb{F}Q_{\Gamma}/I$ be the Brauer configuration algebra associated with a reduced Brauer configuration $\Gamma$ (i.e., truncated vertices $\alpha\in\Gamma_{0}$ occur only in polygons with two vertices). Denote by $\pi:\mathbb{F}Q_{\Gamma}\rightarrow \Lambda$ the canonical surjection then $\pi(x)$ is denoted by $\overline{x}$, for $x\in \mathbb{F}Q_{\Gamma}$.
\end{Nota}

Henceforth, if no confusion arises, we will assume notations $Q$, $I$, and $\Lambda$ instead of $Q_{\Gamma}$, $I_{\Gamma}$, and $\Lambda_{\Gamma}$, for a quiver, an admissible ideal, and the Brauer configuration algebra induced by a fixed Brauer configuration $\Gamma$. \par\bigskip As an example of the application of Algorithm \ref{algoritmo0} consider the following reduced Brauer configuration:

\begin{equation}\label{ex1}
\begin{split}
\Gamma&=(\Gamma_{0},\Gamma_{1},\mu,\mathcal{O}),\\
\Gamma_{0}&=\{0,1\},\\
\Gamma_{1}&=\{V_{1}=\{0,1,1\}, V_{2}=\{0,1\}\},\\
S_{0}&=V_{1}<V_{2},\quad C_{0}=V_{1}<V_{2}<V_{1},\\
S_{1}&=V_{1}<V_{1}<V_{2},\quad C_{1}=V_{1}<V_{1}<V_{2}<V_{1},\\
val(0)&=2,\\
val(1)&=3,\\
\mu(0)&=\mu(1)=1.
\end{split}
\end{equation}

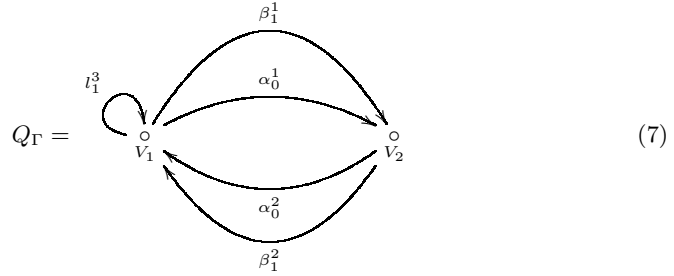
\begin{equation}\label{ex2}
Q_{\Gamma}=\xymatrix@=80pt{\underset{V_{1}}{\circ}\ar@(l,u)@[\blue]^{l_1^3}\ar@/^15pt/@[\blue][r]^{\alpha^{1}_{0}}\ar@/^40pt/@[\blue][r]^{\beta^1_1}&\underset{V_{2}}{\circ}\ar@/^20pt/@[\blue][l]^{\alpha^2_0}\ar@/^40pt/@[\red][l]^{\beta^2_1}}
\end{equation}

The following is a list of special cycles:

\begin{equation}
\begin{split}
C^{1}_{0,V_{1}}&=\alpha^{1}_{0}\alpha^{2}_{0},\\
C^{1}_{0,V_{2}}&=\alpha^{2}_{0}\alpha^{1}_{0},\\
C^{1}_{1,V_{1}}&=l^{3}_{1}\beta^{1}_{1}\beta^{2}_{1},\\
C^{2}_{1, V_{1}}&=\beta^{1}_{1}\beta^{2}_{1}l^{3}_{1},\\
C^{1}_{1,V_{2}}&=\beta^{2}_{1}l^{3}_{1}\beta^{1}_{1}.
\end{split}
\end{equation}

The admissible ideal $I_{\Gamma}$ is generated by the following relations:

\begin{equation}
\begin{split}
C^{k}_{i,V_{j}}&=C^{k'}_{i',V_{j'}}, \quad\text{for all possible values of}\hspace{0.1cm}i, j, k, i', j', k',\\\notag
C^{k}_{i,V_{j}}a&=0,\quad\text{where}\hspace{0.1cm}a\hspace{0.1cm}\text{is the first arrow of the special cycle}\hspace{0.1cm}C^{k}_{i, V_{j}},\\
l^{3}_{1}\alpha^{1}_{0}&=0,\\
\beta^{1}_{1}\alpha^{2}_{0}&=0,\\
\alpha^{2}_{0}l^{3}_{1}&=0,\\
\alpha^{2}_{0}\beta^{1}_{1}&=0,\\
\beta^{2}_{1}\alpha^{1}_{0}&=0.
\end{split}
\end{equation}

The next result regards the dimension of a Brauer configuration algebra.

\addtocounter{prop}{2}
\begin{prop}[\cite{Green}, Proposition 3.13]\label{dimension}
	\textit{Let $\Lambda$ be a Brauer configuration algebra associated with the Brauer configuration $\Gamma$ and let $\mathcal{C}=\{C_{1},\dots, C_{t}\}$ be a full set of equivalence class representatives of special cycles. Assume that for $i=1,\dots,t$, $C_{i}$ is a special $\alpha_{i}$-cycle where $\alpha_{i}$ is a non-truncated vertex in $\Gamma$}. \textit{Then}
	\begin{center}
		$\mathrm{dim}_{\mathbb{F}}\hspace{0.1cm}\Lambda=2|Q_{0}|+\underset{C_{i}\in\mathcal{C}}{\sum}|C_{i}|(n_{i}|C_{i}|-1)$,
		
	\end{center}
	\textit{where $|Q_{0}|$ denotes the number of vertices of $Q$, $|C_{i}|$ denotes the number of arrows in the $\alpha_{i}$-cycle $C_{i}$ and $n_{i}=\mu(\alpha_{i})$}.
\end{prop}

\begin{prop}[\cite{Green}, Proposition 3.6]\label{grading}
	\textit{Let $\Lambda$ be the Brauer configuration algebra associated with a connected Brauer configuration $\Gamma$. The algebra $\Lambda$ has a length grading induced from the path algebra $\mathbb{F}Q$ if and only if there is an $N\in\mathbb{Z}_{>0}$ such that for each non-truncated vertex $\alpha\in\Gamma_{0}$ $val(\alpha)\mu(\alpha)=N$.}
\end{prop}
The following result regards the center of a Brauer configuration algebra.

\addtocounter{teor}{5}
\begin{teor}[\cite{Sierra}, Theorem 4.9]\label{Serra}
	\textit{Let $\Gamma$ be a reduced and connected Brauer configuration and let $Q$ be its induced quiver and let $\Lambda$ be the induced Brauer configuration algebra such that $\mathrm{rad}^{2}\hspace{0.1cm}\Lambda \neq 0$ then the dimension of the center of $\Lambda$ denoted $\mathrm{dim}_{\mathbb{F}}\hspace{0.1cm}Z(\Lambda)$ is given by the formula}:
	
	\begin{equation*}\label{Sierra}
	\begin{split}
	\mathrm{dim}_{\mathbb{F}}\hspace{0.1cm}Z(\Lambda)&=1+\underset{\alpha\in\Gamma_{0}}{\sum}\mu(\alpha)+|\Gamma_{1}|-|\Gamma_{0}|+\#(Loops\hspace{0.1cm} Q)-|\mathscr{C}_{\Gamma}|.
	\end{split}
	\end{equation*}
	
	\textit{where} $|\mathscr{C}_{\Gamma}|=\{\alpha\in\Gamma_{0}\mid val(\alpha)=1, \hspace{0.1cm}and\hspace{0.1cm} \mu(\alpha)>1\}$.
	
\end{teor}

For instance, the dimension of the algebra $\Lambda_{\Gamma}=\mathbb{F}Q_{\Gamma}/I_{\Gamma}$ defined by (\ref{ex1}) and (\ref{ex2}) is given by the following identity:

\begin{equation}
\mathrm{dim}_{\mathbb{F}}\hspace{0.1cm}\Lambda_{\Gamma}=2(2)+2(1)+3(2)=12.\notag
\end{equation}

Note that,

\begin{equation}
\mathrm{dim}_{\mathbb{F}}\hspace{0.1cm}Z(\Lambda_{\Gamma})=1+2+1=4.\notag
\end{equation}

\subsubsection{The Message of a Brauer Configuration}\label{Message}
The notion of labeled Brauer configurations and the message of a Brauer configuration were introduced by Fern\'andez et al in order to define suitable specializations  of some Brauer configuration algebras \cites{Fernandez1, Fernandez}. According to them, since polygons in a Brauer configuration $\Gamma=(\Gamma_{0},\Gamma_{1},\mu,\mathcal{O})$ are multisets, it is possible to assume that any polygon $U\in\Gamma_{1}$ is given by a word $w(U)$ of the form

\begin{equation}\label{word}
\begin{split}
w(U)&=\alpha_{1}^{s_1}\alpha_{2}^{s_2}\dots \alpha^{s_{t-1}}_{t-1}\alpha_{t}^{s_{t}}
\end{split}
\end{equation}

where for each $i$, $1\leq i\leq t$, $s_{i}=\mathrm{occ}(\alpha_{i},U)$.\par\bigskip  

The message is in fact an element of an algebra of words  $\mathscr{W}_{\Gamma}$ associated with a fixed  Brauer configuration such that for a given field $\mathbb{F}$ the word algebra $\mathscr{W}_{\Gamma}$ consists of formal sums of words with the form $\underset{\underset{U\in\Gamma_{1}}{\alpha_{i}\in \mathbb{F}}}{\sum}\alpha_{i}w(U)$, $0w(U)=\varepsilon$ is the empty word, and $1w(U)=w(U)$ for any $U\in\Gamma_{1}$. The product in this case is given by the usual word concatenation. The formal product (or word product)
\begin{equation}\label{message}
\begin{split}
M(\Gamma)&=\underset{U\in\Gamma_{1}}{\prod}w(U)
\end{split}
\end{equation}
is said to be the \textit{message of the Brauer configuration $\Gamma$}.

\par\bigskip

\subsubsection{Mutations}\label{Mutations}
Define a seed $(\Gamma, \mathscr{X})$, where $\mathscr{X}=(x_{1},x_{2},\dots, x_{l})$ with $x_{i}\in\mathbb{F}^{l}$ and a field $\mathbb{F}=\mathbb{Z}_{2}[x]/\langle p(x)\rangle$ for a suitable irreducible polynomial $p(x)$ of degree $n\geq$1.\par\bigskip

$\Gamma=\{\Gamma_{0},\Gamma_{1},\mu,\mathcal{O}\}$ is a Brauer configuration for which;
\begin{equation}\label{mut01}
\begin{split}
\Gamma_{0}&=\{0,1\},\\
\Gamma_{1}&=\{w(0),w(1),\dots,w(l-1)\mid w(i)\in\mathbb{Z}^{n}_{2}\},\\
\mu(0)&=\mu(1)=1.
\end{split}
\end{equation}

For the orientation $\mathcal{O}$ in successor sequences, it is considered the order $w(0)<w(1)<\dots<w(l-1)$.\par\bigskip

A \textit{mutation} $\mathscr{M}(\Gamma,\mathscr{X})=(\Gamma',\mathscr{X}')$ is given by a Brauer configuration $\mathscr{M}(\Gamma)=\Gamma'=(\Gamma'_{0},\Gamma'_{1},\mu',\mathcal{O}')$ and a vector $\mathscr{X}'=(x'_{1},x'_{2},\dots, x'_{l})$ with $x'_{i}\in\mathbb{F}^{l}$, $1\leq i\leq l$ such that:

\begin{equation}\label{mut02}
\begin{split}
\Gamma'_{0}&=\{0,1\},\\
\Gamma'_{1}&=\{w'(0),w'(1),\dots,w'(l-1)\mid w'(i)\in\mathbb{Z}^{n}_{2}\}, \\
\Gamma'_{1}&=\{w(l),w(l+1),\dots, w(2l-1)\mid w(l+i)\in\mathbb{Z}^{n}_{2}\},\hspace{0.1cm}\text{if the seed indices are considered},\\
w'(i)&=\mathscr{M}(x_{i})=w(i-l)+w(i-1),\hspace{0.1cm}\text{if}\hspace{0.1cm}i\not\equiv0\mod{l},\\
w'(i)&=w(i-l)+\mathscr{H}(w(i-1)),\hspace{0.1cm}\text{if}\hspace{0.1cm}i\equiv0\mod{l},\\
\mu(0)&=\mu(1)=1.
\end{split}
\end{equation}
Where $\mathscr{H}(x_{1},x_{2},\dots, x_{l})=(\mathcal{T}(x_{1}),\mathcal{T}(x_{2}),\dots,\mathcal{T}(x_{l}))$ with $\mathcal{T}(x_{j})=\lambda_{j}x_{j}+v_{j,0}$, for some suitable, $\lambda_{j},v_{j,0}\in\mathbb{F}$, we let $v_{0}$ denote the shift vector such that $v_{0}=(v_{1,0},\dots, v_{l,0})$.\par\bigskip

For successor sequences, the orientation $\mathcal{O'}$ is defined in such a way that,
\begin{equation}
 w'(0)<w'(1)<\dots<w'(l-1)\notag
 \end{equation}

It turns out that according to the indices assumed for the original seed, the $i$th mutation $\mathscr{M}^{i}$ has the form:

\begin{equation}
\mathscr{M}^{i}=\{w(i(l-1)+1),w(i(l-1)+2),\dots, w(i(l-1)+l)\mid w(i)\in\mathbb{Z}^{n}_{2}\}.\notag
\end{equation} 

Brauer configurations obtained from mutations are said to be \textit{Brauer clusters}. Polygons are called \textit{cluster polygons}.\par\bigskip

For a fixed positive integer $m_{0}$, the $m_{0}$-\textit{Brauer cluster}, $\Phi^{m_{0}}$ is a Brauer configuration, 

\begin{equation}\label{mutation}
\Phi^{m_{0}}=(\Phi^{m_{0}}_{0},\Phi^{m_{0}}_{1},\mu,\mathcal{O})
\end{equation}

 such that:

\begin{enumerate}
\item $\Phi^{m_{0}}_{0}=\{0,1\}$,
\item $\Phi^{m_{0}}_{1}$ consists of messages of the Brauer clusters $\mathscr{M}^{i}$, $i\geq1$,
\item $\mu(0)=\mu(1)=1$,
\item If $M(i)$ denotes the message associated with the $i$th Brauer cluster then for successor sequences, it is assumed the order
$M(0)<M(1)<\dots<M(m_{0})$.

\end{enumerate}

As an example, in the sequel we describe mutations of the following seed $(\Gamma,\mathscr{X}=\{1,x\})$, where $\mathbb{F}=\mathbb{Z}_{2}[x]/\langle x^{2}+x+1\rangle$.
\begin{equation}
\begin{split}
\Gamma&=(\Gamma_{0},\Gamma_{1},\mu,\mathcal{O}),\\
\Gamma_{0}&=\{0,1\},\\
\Gamma_{1}&=\{1=01, x=10\},\\
\mu(0)&=\mu(1)=1,\\
v_{1,0}&=1,\quad \lambda_{1}=x,
\end{split}
\end{equation}
\par\bigskip
Successor sequences at $0$ and $1$ are defined as follows:

\begin{equation}
0:1<x,\quad
1:1<x.\notag
\end{equation}
\par\bigskip
The first mutation $\mathscr{M}(\Gamma)=(\Gamma',\mathscr{X}'=\{x, x+1\})$ is defined as follows:

\begin{equation}
\begin{split}
\Gamma'&=(\Gamma'_{0},\Gamma'_{1},\mu',\mathcal{O}'),\\
\Gamma'_{0}&=\{0,1\},\\
\Gamma'_{1}&=\{x=10, x+1=11\},\\
\mu'(0)&=\mu'(1)=1,\\
v_{2,0}&=1,\quad \lambda_{2}=x.
\end{split}
\end{equation}

In this case 1 is the only nontruncated vertex with a successor sequence of the form:

\begin{equation}
1:x<x+1<x+1.\notag
\end{equation}
\par\bigskip

The second mutation $\mathscr{M}^{2}=(\Gamma'',\mathscr{X}''=\{x,1\})$ is described as follows:
\begin{equation}
\begin{split}
\Gamma''&=(\Gamma''_{0},\Gamma''_{1},\mu'',\mathcal{O}''),\\
\Gamma''_{0}&=\{0,1\},\\
\Gamma''_{1}&=\{x=10, 1=01\},\\
\mu''(0)&=\mu''(1)=1.
\end{split}
\end{equation}

In this case, the successor sequences  at 0 and 1 are
\begin{equation}
0: x<1,\quad 1:x<1.\notag
\end{equation}

The 2-Brauer cluster $\Phi^{2}=(\Phi^{2}_{0},\Phi^{2}_{1},\mu,\mathcal{O})$ has the following form:

\begin{equation}
\begin{split}
\Phi^{2}_{0}&=\{0,1\},\\
\Phi^{2}_{1}&=\{M(0)=0110,M(1)=1011,M(2)=1001\},\\
\mu(0)&=\mu(1)=1.
\end{split}
\end{equation}

For successor sequences, it holds that:\par\bigskip
\begin{equation}
M(0)<M(1)<M(2).
\end{equation}
\par\bigskip
The following $Q_{\Gamma}, Q_{\Gamma'}, Q_{\Gamma''}$ and, $Q_{\Phi^{2}}$ are the Brauer quivers of all these mutations and of the 2-Brauer cluster $\Phi^{2}$.
\begin{equation}\label{mut0}
Q_{\Gamma}=\xymatrix@=80pt{\underset{1}{\circ}\ar@/^15pt/@[\blue][r]^{\alpha^{1}_{0}}\ar@/^40pt/@[\blue][r]^{\alpha^2_0}&\underset{x}{\circ}\ar@/^20pt/@[\blue][l]^{\beta^1_1}\ar@/^40pt/@[\red][l]^{\beta^2_1}}\notag
\end{equation}

\begin{equation}\label{mut1}
Q_{\Gamma'}=\xymatrix@=80pt{\underset{x}{\circ}\ar@/^15pt/@[\blue][r]^{\beta^{1}_{1}}&\underset{x+1}{\circ}\ar@(dr,ru)@[\green]_{l_1^1}\ar@/^15pt/@[\red][l]^{\beta^2_1}}
\end{equation}

\begin{equation}\label{mut2}
Q_{\Gamma''}=\xymatrix@=80pt{\underset{x}{\circ}\ar@/^15pt/@[\blue][r]^{\alpha^{1}_{0}}\ar@/^40pt/@[\blue][r]^{\alpha^2_0}&\underset{1}{\circ}\ar@/^20pt/@[\blue][l]^{\beta^1_1}\ar@/^40pt/@[\red][l]^{\beta^2_1}}\notag
\end{equation}

\begin{equation}\label{mutfinal0}
Q_{\Phi^{2}}=\xymatrix@=80pt{*+[o][F-]{M_{0}}\ar@/^0pt/@[\blue][r]^{\beta^{1}_{0}}\ar@/^40pt/@[\blue][r]^{\alpha^1_1}\ar@(u, l)@[\blue]_{l_1^1}\ar@(l, d)@[\blue]_{l_0^1}&*+[o][F-]{M_{1}}\ar@(d, r)@[\blue]_{l_0^2}\ar@(l, d)@[\blue]_{l_1^2}\ar@/^0pt/@[\blue][r]^{\beta^2_0}\ar@/^40pt/@[\red][r]^{\alpha^2_1}&*+[o][F-]{M_{2}}\ar@(u, r)@[\blue]^{l_1^3}\ar@/^40pt/@[\red][ll]^{\alpha^3_1}\ar@(d, r)@[\blue]_{l_0^3}\ar@/^60pt/@[\red][ll]^{\beta^3_0}\ar@(d, r)@[\blue]_{l_0^3}}
\end{equation}

The admissible ideal $I$ in the Brauer configuration algebra $\Lambda_{(\Phi^{n},\mathscr{X})}$ is generated by the following relations:

\begin{enumerate}
\item $(l^{i}_{j})^{2}$,\quad $l^{i}_{0}l^{j}_{1}$, for all possible values of $i$ and $j$,
\item $\alpha^{i}_{1}\beta^{r}_{0}$,\quad$\alpha^{i}_{1}l^{r}_{0}$,\quad$l^{r}_{1}\alpha^{i}_{1}$, for all possible values of $i$ and $r$,
\item $\beta^{r}_{0}\alpha^{i}_{1}$,\quad$\beta^{r}_{0}l^{s}_{1}$,\quad$l^{s}_{j}\beta^{r}_{0}$, for all possible values of $i, r$, and $s$.

\end{enumerate}

\begin{teor}\label{MUT}
For $m_{0}>1$ the Brauer configuration algebra $\Lambda_{\Phi^{m_{0}}}$ induced by the Brauer configuration $\Phi^{m_{0}}$ (see $\mathrm{(\ref{mutation})}$) is connected and reduced.
\end{teor}

\textbf{Proof.} Since each polygon $M(i)\in\Phi^{m_{0}}_{1}$ contains at least one 0 and at least one 1, we conclude that $\Phi^{m_{0}}$ does not contain truncated vertices. The result follows provided that $\underset{i=0}{\overset{m_{0}}{\cap}}M(i)\neq\varnothing$.\hspace{0.5cm}$\square$

\begin{teor}\label{mutteor}
If $\mathbb{F}$ is a finite field then the set $M_{\Phi}$ of $m_{0}$-Brauer clusters obtained by mutation is finite, and the special cycle of maximal length has one of the two shapes $\mathrm{(\ref{mut1s})}$ or $\mathrm{(\ref{mut1sa})}$. Moreover, if $\Phi^{i}\neq\Phi^{j}$ whenever $i\neq j$ then there exists $m\in\mathbb{N}$ such that  $(\Gamma, \mathscr{X})= (\Gamma^{m+1}, \mathscr{X}^{m+1})$ for a given initial seed  $(\Gamma, \mathscr{X})$.

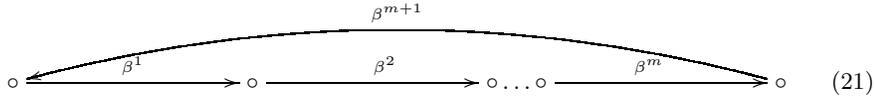
\begin{equation}\label{mut1ss}
\xymatrix@=80pt{{\circ}\ar@/^0pt/@[\blue][r]^{\beta^{1}}&{\circ}\ar@/^0pt/@[\red][r]^{\beta^2}& \circ\dots \ar@/^0pt/@[\red][r]^{\beta^m}\circ&\circ \ar@/_20pt/@[\red][lll]_{\beta^{m+1}}}
\end{equation}

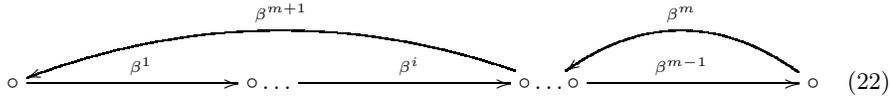
\begin{equation}\label{mut1saa}
\xymatrix@=80pt{{\circ}\ar@/^0pt/@[\blue][r]^{\beta^{1}}&{\circ}\dots\ar@/^0pt/@[\red][r]^{\beta^i}& \circ\dots \ar@/^0pt/@[\red][r]^{\beta^{m-1}}\ar@/_20pt/@[\red][ll]_{\beta^{m+1}}\circ&\circ \ar@/_20pt/@[\red][l]_{\beta^{m}}}
\end{equation}

where $\beta^{m}\beta^{m+1}\neq0$.

\end{teor}

\textbf{Proof.} Every mutation $(\Gamma^{i}, \mathscr{X}^{i})$ is uniquely determined by $\mathscr{X}^{i}$. Since, $\mathscr{X}^{i}\in \mathbb{F}^{l}$ for any $i\in\mathbb{N}$ then $|M_{\Phi}|\leq |\mathbb{F}^{l}|$. We are done.\hspace{0.5cm}$\square$

\addtocounter{Nota}{5}
\begin{Nota}\label{probability}
\normalfont
For $m_{0}>1$ and a fixed seed, the dimension of an algebra $\Lambda_{\Phi^{m_{0}}}$ and its center $Z(\Lambda_{\Phi^{m_{0}}})$ can be estimated by using some statistical methods. For instance, we use a sample of $10^{6}$ random seeds in order to obtain confidence intervals for these values. Such samples allow us to infer that if $m_{0}=10$ then
\begin{equation}
P_{r}(7499\leq\mathrm{dim}_{\mathbb{F}}\hspace{0.1cm}\Lambda_{\Phi^{m_{0}}}\leq 8067)\geq 0,99\quad\text{and}\quad P_{r}(199\leq\mathrm{dim}_{\mathbb{F}}\hspace{0.1cm}Z(\Lambda_{\Phi^{m_{0}}})\leq 221)\geq 0,99.\notag
\end{equation}

Where $P_{r}(X)$ denotes the probability of an event $X$.
\end{Nota}

\subsection{Deterministic and Non-Deterministic  Automata}

In this section, we recall definitions of deterministic and non-deterministic automata as Rutten et al present in \cite{Ballester}. We use these definitions to interpret Brauer configuration algebras as automata with acceptance language given by relations generating suitable admissible ideals.\par\bigskip

Given an alphabet $A$, a \textit{deterministic automaton} is a pair $(X,\alpha)$ consisting of a (possibly infinite) set $X$ of states and a transition function $\alpha: X\rightarrow X^{A}$. The following is an illustration of this kind of transitions, where $\alpha(x)(a)=y=x_{a}$ \cite{Ballester}.

\begin{equation}\label{fexample}
\xymatrix@=80pt{*+[o][F-]{x}\ar@/^0pt/@[\blue][r]^{a}&*+[o][F-]{y}}\notag
\end{equation}
\par\bigskip
If $\varepsilon$ denotes the empty word then $x_{\varepsilon}=x$, for any $x\in X$ and $x_{wa}=\alpha(x_{w})(a)$ with $w\in A^{*}$.\par\bigskip

A deterministic automaton can be decorated by means of a coloring function $c:X\rightarrow 2=\{0,1\}$ such that $c(x)=1$ if $x$ is an \textit{accepting (or final)} state, $c(x)=0$, if $x$ is a \textit{non-accepting} state. A triple $(X, c,\alpha)$ is said to be a deterministic
colored automaton. In the following diagram an accepting state is denoted with a double circle, $x$ is an accepting state, whereas $y$ is a non-accepting state.

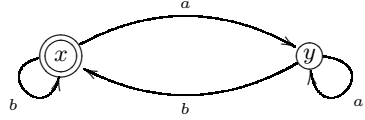
\begin{equation}\label{fexamplee}
\xymatrix@=80pt{*++[o][F=]{x}\ar@/^15pt/@[\blue][r]^{a}\ar@(l, d)@[\blue]_{b}&*+[o][F-]{y}\ar@/^15pt/@[\red][l]^{b}\ar@(r, d)@[\blue]^{a}}\notag
\end{equation}
\par\bigskip

Given a deterministic colored automaton $(X, c,\alpha)$ and a state $x\in X$, the set 

\begin{equation}O_{c}(x)=\{w\in A^{*}\mid c(x_{w})=1\}\notag
\end{equation}
\par\bigskip
is called the \textit{language accepted} or recognized by the automaton $(X, c,\alpha)$ starting from the state $x$. A deterministic automaton can also has an \textit{initial state} $x\in X$ represented by a function $x:1=\{0\}\rightarrow X$. The triple $(X, x,\alpha)$ is said to be a \textit{deterministic pointed automaton}.\par\bigskip

A \textit{non-deterministic automaton} is a pair $(X,\alpha)$ consisting of a set $X$ (possibly infinite) of states and a transition function $\alpha: X\rightarrow P_{w}(X)^{A}$, that assigns to each letter and to each state a finite set of states \cite{Ballester}. If to each state it is assigned a single new state, the definition of deterministic automaton is recovered. As in the deterministic case, a state $x$ in a non-deterministic automaton can be either accepting ($c(x)=1$) or non-accepting ($c(x)=0$). And $x_{\varepsilon}=\{x\}$, $x_{w_{a}}=\bigcup\{y_{a}\mid y\in x_{w}\}$. A triple $(X, c,\alpha)$ is called a \textit{colored non-deterministic automaton}.\par\bigskip

\par\bigskip

\begin{centering}
\textbf{A Regular Language associated with a Brauer Configuration Algebra}\par\bigskip
\end{centering}

Automata associated with path algebras have been studied by Rees \cite{Rees}, who introduced some automata associated with some string algebras, such automata were used by her in order to describe indecomposable representations over these types of algebras, she points out that the set of strings defining representations of string algebras and many other bounded path algebras constitute a regular set. In this section, we follow the ideas of Rees to describe an automaton associated with a Brauer configuration algebra. \par\bigskip

Values of the map $c: Q_{1}\rightarrow 2$ can be obtained by endowing to the successor sequences a length-lexicographic order, in this case, both $\Gamma_{0}$ and $\Gamma_{1}$ are well-ordered sets with partial orders $\prec$ and $<$ respectively. In such a way that initial states in the corresponding automaton are given by minimal successor sequences (see Algorithm \ref{algoritmo0}). Note that, if  $S_{a,U}$ denotes a successor sequence starting in a polygon $U$ with $a\in U$ and if $|S_{\alpha,V_{1}}|=|S_{\alpha',V_{1}}|=|S_{\alpha,V_{2}}|=|S_{\alpha',V_{2}}|$, $\alpha,\alpha'\in V_{1}\cap V_{2}$ and $\alpha\prec\alpha'$ then $c(S_{\alpha, V_{1}})=1$ and $c(S_{\alpha', V_{1}})=c(S_{\alpha, V_{2}})=c(S_{\alpha', V_{2}})=0$.\par\bigskip

A  Brauer configuration algebra $\Lambda_{\Gamma}$ induced by a Brauer configuration $\Gamma=(\Gamma_{0},\Gamma_{1},\mu,\mathcal{O})$ has associated a regular language $L_{\Gamma}=A^{*}_{\Gamma}/\sim$, where the alphabet $A_{\Gamma}=\{x^{i}_{\alpha}\mid \alpha\in\Gamma_{0}, 1\leq i\leq val(\alpha)\}$, each letter $x^{i}_{\alpha}$ corresponds to a unique arrow in $(Q_{\Gamma})_{1}$. And each path $P\in Q_{\Gamma}$ corresponds to a word $w\in L_{\Gamma}$. \par\bigskip

Two words $w, w'\in L_{\Gamma}$ are equivalent (i.e., $w\sim w'$) if their corresponding paths are equivalent as elements of the Brauer configuration algebra. In this case, if $S_{\alpha}$ is a successor sequence associated with the vertex $\alpha\in \Gamma_{0}$, then $w_{S_{\alpha}}$ denotes the word associated with the corresponding special cycle up to equivalence. If $\mathrm{min}\hspace{0.1cm}S_{\alpha}=U\in\Gamma_{1}$ then  $w_{S_{\alpha}}\in O_{c}(U)$ (final vertices of special cycles are final states up to equivalence).

\par\bigskip
In the associated automaton of a Brauer configuration algebra, polygons are states. Actually, all states we represent are accepted states. Transition between states is given by the order $<$, in other words, if $x^{\alpha}_{i}$ is the letter associated with an arrow $U_{i}<U_{i+1}$, that is, $U_{i}\stackrel{x^{i}_{\alpha}}{\longrightarrow}U_{i+1}\in (Q_{\Gamma})_{1}$ then $x^{i}_{\alpha}$ is a transition from $U_{i}$ to $U_{i+1}$. Note that, if ${x^{i}_{\alpha}}{x^{j}_{\alpha'}}$ belongs to the admissible ideal $I$, with $\Lambda_{\Gamma}=kQ_{\Gamma}/I$ then it is not accepted as word in $L_{\Gamma}$ if $\alpha\neq \alpha'$. \par\bigskip

Note that, according to the automaton associated with the Brauer configuration defined by (\ref{ex1}) and (\ref{ex2}), it holds that $V_{1}$ is the initial and final state and

\begin{equation}
\begin{split}
c(C^{k}_{i,V_{j}})&=0,\quad\text{if}\hspace{0.1cm}, i\neq0, j\neq1, k\neq1. \\\notag
c(C^{k}_{i,V_{j}}a)&=0,\quad\text{if}\hspace{0.1cm} a\hspace{0.1cm}\text{is the first arrow of}\hspace{0.1cm}C^{k}_{i,V_{j}}\hspace{0.1cm}\text{for all possible values of}\hspace{0.1cm}i,j\hspace{0.1cm}\text{and}\hspace{0.1cm}k.\\
c(\alpha^{i}_{j}\beta^{i'}_{j'})&=c(\beta^{i}_{j}\alpha^{i'}_{j'})=0, \text{for all the possible values of}\hspace{0.1cm}i, j, i', j'.
\end{split}
\end{equation}

\subsection{Enumeration of Gelfand-Tsetlin Patterns}

In this section, we recall some well known results regarding enumeration of GT-patterns. And introduce the notion of heart of a GT-pattern which allows us to define posets and marked order polytopes associated with the number of some GT-patterns \cites{DeLoera, Fisher1, Fisher2, Futorny, RamirezF}.\par\bigskip

If $\lambda=(\lambda_{(n,1)},\lambda_{(n,2)},\lambda_{(n,3)},\dots, \lambda_{(n,n)})$ is an integer partition and $V(\lambda)$ is a finite dimensional irreducible representation of $\mathrm{gl}_{n}$ with highest weight $\lambda$ then a basis of $V(\lambda)$ is parametrized by GT-patterns $\mathscr{T}=\mathscr{T}((n,1),(n,n),(1,1))$ associated with $\lambda$. These are arrays of integer row vectors with the shape:

\par\bigskip

\begin{centering}

$\mathscr{T}=
  \begin{array}{ccccccc}
    \lambda_{(n,1)} & \lambda_{(n,2)} & \lambda_{(n,3)} &\dots & \dots &\dots& \lambda_{(n,n)}   \\
     & \lambda_{(n-1,1)} & &\lambda_{(n-1,2)}&\dots&\lambda_{(n-1,n-1)}&\\
     &\ddots&\dots & \dots&\dots & \iddots& \\
     &&&&&&\\
     &&&\lambda_{(1,1)} &&&
        
  \end{array}$ \par\bigskip

  \end{centering}

such that the upper row coincides with $\lambda$ and the following conditions hold:

\begin{equation}\label{GTconditions}
\lambda_{(k,i)}\geq\lambda_{(k-1,j)},\quad\lambda_{(k-1,i)}\geq\lambda_{(k,i+1)}.
\end{equation}
\par\bigskip
This setting together with terms of the form $l_{(k, i)}=\lambda_{(k ,i)}-i+1$ allow establishing the existence of a basis of $V(\lambda)$ parametrized by $\mathscr{T}$ and such that the action of generators of $\mathrm{gl}_{n}$ is given  by some Gelfand-Tsetlin formulas \cites{Futorny, RamirezF}.
\par\bigskip
A GT-pattern with first row of the form $1, 2,\dots, n$ is called a monotone triangle of length $n$. It is also known that there is a bijection
between GT-patterns with first row $\lambda_{n}\leq \lambda_{n-1}\leq\dots\leq\lambda_{1}$ and column-strict plane partitions of type $\lambda=(\lambda_{1},\lambda_{2},\dots,\lambda_{n})$ ($\lambda_{i}$ parts in row $i$) and largest part $\leq n$ \cites{Fisher1, Fisher2}. \par\bigskip

The following theorem regarding monotone triangles was given by Zeilberger in 1996.
\addtocounter{teor}{1}
\begin{teor}[\cite{Zeilberger}, Main Theorem]
The number of monotone triangles of length $n$ with bottom entry $a_{(n, n)}=r$ is equal to
\begin{equation}
\binom{2n-2}{n-1}\binom{n+r-2}{n-1}\binom{2n-r}{n-1}A_{n-1}, \quad\text{with}\quad A_{n}=\underset{i=0}{\overset{n-1}{\prod}}\frac{(3i+1)!}{(n+1)!}.\notag
\end{equation}
\end{teor}
\par\bigskip
Note that, $A_{n}(x)=\underset{\mathscr{T}}{\sum}x^{s(\mathscr{T})}$, where $s(\mathscr{T})$ denotes the number of standard elements of $\mathscr{T}$ such that $t_{(i-1,j-1)}<t_{(i,j)}<t_{(i-1,j)}$, for $2\leq i\leq j\leq n$. In particular, $A_{n}(2)=2^{\binom{n}{2}}$.
\section{Gelfand-Tsetlin equations and Diophantine Equations of Type $\mathcal{D}(n_{1},n_{2},\mathcal{K}_{m})$}
In this section, Gelfand-Tsetlin equations and Gelfand-Tsetlin numbers are introduced. Solutions of this kind of equations are proposed via some order marked polytopes. Whereas, solutions of diophantine equations of type  $\mathcal{D}(n_{1},n_{2},\mathcal{K}_{m})$ are obtained by using some suitable non-deterministic automata associated with mutations of Brauer configurations.
\subsection{The Heart of a GT-Pattern}

For $n\geq4$, a subarray $\mathscr{T}((n-1,1),(n-1,n-2),(2,1))$ of $\mathscr{T}((n,1),(n, n),(1,1))$ is said to be the \textit{heart} or \textit{GT-heart} of $\mathscr{T}((n,1),(n, n),(1,1))$ denoted $\mathcal{H}_{n}(\mathscr{T})$.\par\bigskip

For  $n\geq 4$ fixed, it is possible to define an order $\unlhd$ on the set of GT-hearts over $\mathrm{gl}_{n}$ (associated with GT-patterns with fixed first row) whose covers are defined as follows:
\par\bigskip

$\mathcal{H}_n(\mathscr{T})\unlhd\mathcal{H}_{n}(\mathscr{T}')$ if $\mathcal{H}_{n}(\Lambda)\backslash\{t_{(i,j)}\}=\mathcal{H}_{n}(\Lambda)\backslash\{t'_{(i,j)}\})$ and $t'_{(i,j)}=t_{(i,j)}+a$, with $a\in\{0,1\}$, for some $i,j$.

\par\bigskip
 If $T(\mathcal{H}_{n})$ denotes the set of Gelfand-Tsetlin arrays with $\mathcal{H}_{n}$ as heart, then
two hearts $\mathcal{H}_{n}$ and $\mathcal{H}'_{n}$ are said to be equivalent denoted $\mathcal{H}_{n}\cong\mathcal{H}'_{n}$ provided that $|T(\mathcal{H}_{n})|=|T(\mathcal{H}'_{n})|$. Thus, $\mathscr{H}_{(n, r)}$ is a poset endowed with an equivalence relation.\par\bigskip

We let $(\mathscr{H}_{(n, r)},\unlhd)$ or simply $\mathscr{H}_{(n, r)}$ denote the poset of hearts associated with a given Gelfand-Tsetlin pattern with a weight $\lambda=(\lambda_{(n,1)},\dots,\lambda_{(n, n)})$ of length $n$ and such that $|\lambda_{(n, i)}-\lambda_{(n, i-1)}|=r$. \par\bigskip

For instance, if $\lambda=(\lambda_{(4,1)}=m,\lambda_{(4,2)}=m-3,\lambda_{(4,3)}=m-6,\lambda_{(4,4)}=m-9)$ is a suitable weight for $n=4$, then the following is a three-point chain of $\mathscr{H}_{(4,3)}$: 

\par\bigskip
\begin{centering}

$
  \begin{array}{ccc}
   
     2& -1&\\
     &2&        
  \end{array}\unlhd \begin{array}{ccc}
   
     3& -1&\\
     &2&        
  \end{array}\unlhd \begin{array}{ccc}
   
     4& -1&\\
      &2&      
  \end{array}$ \par\bigskip

  \end{centering}

We have the following result for GT patterns with first row of the form: 

\begin{equation}
w=(n,n-r,n-2r,\dots, n-r(n-1)).\notag
\end{equation}

\begin{teor}\label{NGT}
The number $g_{(n, r)}$ of GT-patterns over $\mathrm{gl}_{n}$ defined by a weight vector of the form
$w=(n, n-r, n-2r,\dots, n-r(n-1))$ is given by the formula $g_{(n, r)}=(r+1)^{\binom{n}{2}}$.

\end{teor}

\textbf{Proof.} It suffices to observe that for $n>2$, it holds that $g_{(n, r)}=\frac{g^{2}_{(n-1,r)}(r+1)}{g_{(n-2,r)}}$. \hspace{0.5cm}$\square$

\par\bigskip
If $|\mathcal{H}_{n}(\mathscr{T})|$ denotes the number of Gelfand-Tsetlin patterns $\mathscr{T}((n,1),(n, n),(1,1))$ with the same heart $\mathcal{H}_{n}$ then the following result holds.
\addtocounter{prop}{6}
\begin{prop}\label{NGT1}
\textit{If $n=4$, $r>1$, and $\mathcal{H}_{n},\mathcal{H}'_{n}\in\mathscr{H}_{(n, r)}$ are such that $\lambda_{(n-1,2)}=\lambda'_{(n-1,2)}$ and
$\lambda_{(n-2,1)}=\lambda'_{(n-2,1)}$ then $|T(\mathcal{H}_{n})|=|T(\mathcal{H}'_{n})|$.}

\end{prop}

\textbf{Proof.} If the different Gelfand-Tsetlin patterns are obtained by keeping without changes the heart and reducing in one unit only one entry $\lambda_{(j, j)}$ then the number of such Gelfand-Tsetlin arrays is given by adding some suitable integers in the form:
\begin{equation}\label{sequence}
S_{gt}(n_{1},r)=n_{1}q_{1}+q_{2}+\dots+q_{(2r)}+q_{(2r+1)}=n_{1}q_{1}+\underset{j=2}{\overset{2r+1}{\sum}}q_{j}
\end{equation}
in particular, $q_{(2r)}=3$ and $q_{(2r+1)}=1$. In the case, $\lambda_{(2,2)}=\lambda'_{(2,2)}=r-1=\mathrm{max}(\lambda_{(2,2)})$, $1\leq n_{1}\leq r+1$, and
\begin{equation}\label{sequence1}
\begin{split}
q_{1}&=[(2r+1)+r(r+1)+t_{(r-1)}],\\
q_{2}&=[2r+r^{2}+t_{(r-1)}],\\
q_{3}&=[(2r-1)+r(r-1)+t_{(r-1)}],\\
\vdots&=\vdots,\\
q_{r}&=[(2r+1)-r+1+2r+t_{(r-1)}]=[3r+t_{(r-1)}+2],\\
q_{(r+1)}&=t_{(r+1)},\\
q_{(r+2)}&=t_{r},\\
\vdots&=\vdots,\\
q_{2r}&=t_{2},\\
q_{(2r+1)}&=t_{1}=1.
\end{split}
\end{equation}

For arrays with $\lambda_{(2,2)}=r-2$ then $1\leq n_{1}\leq n_{1}+1$, $q_{1}=[2r+r(r)+t_{(r-1)}]$, $q_{2}=[(2r-1)+r(r-1)+t_{(r-1)}]$ and the construction of the remain $q_{i}$ goes on until reaching all arrays for $-1\leq \lambda_{(2,2)}\leq r-1$. If $\lambda_{(2,2)}=-1$ then
we have $2r+1$ sums of the form $n_{1}t_{(r+1)}+\underset{j=1}{\overset{r}{\sum}}t_{r}$, $1\leq n_{1}\leq 2r+1$, $n_{1}=2r+1$ if in the heart $\lambda_{(3,1)}=2r$, $n_{2}=2r$, if in the heart $\lambda_{(3,1)}=2r-1$ and so on until the case for which $\lambda_{(3,1)}=0$ and $n_{1}=1$.\hspace{0.5cm}$\square$

\subsubsection{Marked Posets and Marked Polytopes}

In this section, we describe a special class of marked posets and marked polytopes  introduced by Fourier in \cite{Fourier}. We also prove that some posets of type $\mathscr{H}_{(n, r)}$ are marked by the number of some suitable Gelfand-Tsetlin arrays. \par\bigskip

Let $(\mathscr{P},\preceq)$ be a finite poset and a subset $A$ of $\mathscr{P}$ containing at least all maximal and all minimal elements 
of $\mathscr{P}$, we set:

\begin{equation}\label{markedposet}
Q_{A}=\{\lambda\in\mathbb{Z}^{A}_{\geq0}\mid \lambda_{a}\leq \lambda_{b}\hspace{0.1cm}\text{if}\hspace{0.1cm} a\preceq b\}
\end{equation} 

 and call the triple $(\mathscr{P},A,\lambda)$ a \textit{marked poset}. Then the \textit{marked chain polytope} associated with $\lambda\in Q_{A}$ is defined:
 
 \begin{equation}\label{chainpolytope}
 \mathcal{C}(\mathscr{P},A)=\{s\in\mathbb{R}^{\mathscr{P}\backslash A}_{\geq0}\mid s_{x_{1}}+s_{x_{2}}+\dots+s_{x_{n}}\leq \lambda_{b}-\lambda_{a}\hspace{0.1cm}\text{for all chains}\hspace{0.1cm}a\preceq x_{1}\preceq\dots\preceq x_{n}\preceq b\}\notag
 \end{equation}

while the marked order polytope is defined as

\begin{equation}\label{markedorderpolytope}
\mathcal{O}(\mathscr{P},A)_{\lambda}=\{s\in\mathbb{R}^{\mathscr{P}\backslash A}_{\geq0}\mid s_{x}\leq s_{y},\lambda_{a}\leq s_{x}\leq \lambda_{b},\hspace{0.1cm}\text{for all}\hspace{0.1cm}a\preceq x\preceq b, x\preceq y\}\notag
\end{equation}

where $a, b\in A$, $x_{i},x, y\in\mathscr{P}\backslash A$ (the definition of marked order polytope is also valid in case that order $\leq$ in $\mathbb{R}$ reversed the order $\preceq$ in $\mathscr{P}$).

\par\bigskip

The marked poset $(\mathscr{P},A,\lambda)$ is \textit{regular} if for all $a\neq b$ in $A$, $\lambda_{a}\neq \lambda_{b}$ and there are no obviously redundant relations.\par\bigskip

Let $(\mathscr{P},A,\lambda)$ be a regular marked poset then the number of facets in the marked order polytope is equal to the number of cover relations in $\mathscr{P}$.\par\bigskip

If $c(a, b)$ denotes the number of saturated chains $a\preceq x_{1}\preceq x_{2}\preceq\dots\preceq x_{p}\preceq b$, $x_{i}\in\mathscr{P}\backslash A$, then the number of facets in the marked chain polytope is equal to $|\mathscr{P}\backslash A|$+$\underset{a\preceq b}{\sum}c(a, b)$.

\par\bigskip

In the $\mathrm{sl}_{n}$ case, for each $n$ and $\lambda=(\lambda_{1}\geq\dots\geq \lambda_{n})$, let $(\mathscr{P},A,\lambda)$ be the following marked poset $\mathscr{P}_{\mathrm{sl}_{n}}=\{x_{(i,j)}\mid 0\leq i\leq n, 1\leq j\leq n+1-i\}$ with cover relations $x_{(i-1,j+1)}\geq x_{(i,j)}\geq x_{(i-1,j)}$ \cite{Fourier}.

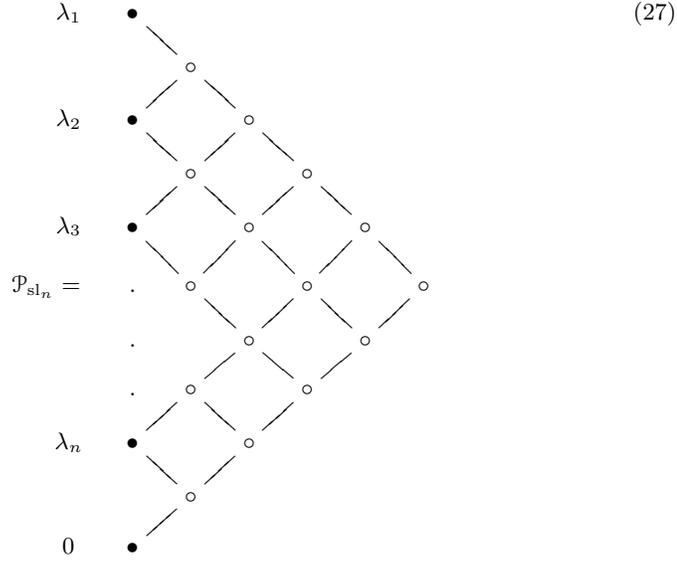
\begin{equation}
{\xymatrix@=14pt@R=3mm@C=4mm{
\lambda_{1}\hspace{-6mm}&{\bullet}\ar@{-}[rd]& &&&&\\
&& {\circ}\ar@{-}[rd]&&&&\\
\lambda_{2}\hspace{-6mm}&{\bullet}\ar@{-}[rd]\ar@{-}[ru]&&{\circ}\ar@{-}[rd]&&&\\
&& {\circ}\ar@{-}[rd]\ar@{-}[ru]&&{\circ}\ar@{-}[rd]&&\\
\lambda_{3}\hspace{-6mm}&{\bullet}\ar@{-}[rd]\ar@{-}[ru]&&{\circ}\ar@{-}[rd]\ar@{-}[ru]&&{\circ}\ar@{-}[rd]&\\
\mathscr{P}_{\mathrm{sl}_{n}}=&.&{\circ}\ar@{-}[rd]\ar@{-}[ru]&&{\circ}\ar@{-}[rd]\ar@{-}[ru]&&{\circ}\\
&.&&{\circ}\ar@{-}[rd]\ar@{-}[ru]&&{\circ}\ar@{-}[ru]&\\
&.&{\circ}\ar@{-}[rd]\ar@{-}[ru]&&{\circ}\ar@{-}[ru]&&\\
\lambda_{n}\hspace{-6mm}&{\bullet}\ar@{-}[rd]\ar@{-}[ru]&&{\circ}\ar@{-}[ru]&&&\\
&& {\circ}\ar@{-}[ru]&&&&\\
0\hspace{-6mm}&{\bullet}\ar@{-}[ru]&&&&&
}}
\end{equation}
\par\bigskip
Let $\lambda=(\lambda_{1}\geq\lambda_{2}\geq\dots\geq\lambda_{n}\geq\lambda_{n+1}=0)$ and let $x_{(0,j)}$ be marked with $\lambda_{j}$, then:

\begin{equation}
\begin{split}
|\{\text{facets in the marked order polytope}\}|&=n(n+1),\\
|\{\text{facets in the marked chain order polytope}\}|&=\frac{n(n-1)}{2}+\underset{i=1}{\overset{n}{\sum}}iC_{n-i},\\
\end{split}
\end{equation}
where $C_{n-i}$ is the corresponding Catalan number.\par\bigskip

The associated marked order polytope is known as the Gelfand-Tsetlin polytope associated with the partition $\lambda$.
\par\bigskip

The following result regards posets of type $(\mathscr{H}_{(n, r)},\unlhd)$.

\begin{prop}\label{NGT2}
\textit{For $n=4$, and $r>1$, if $\mathcal{Q}$ denotes the set of numbers $q_{i}$ given by identities $\mathrm{(\ref{sequence})}$ and $\mathrm{(\ref{sequence1})}$. Then $\mathcal{Q}\subset Q_{\mathscr{H}_{(n,r)}}$ (see $\mathrm{(\ref{markedposet})}$).}

\end{prop}

\textbf{Proof.} Consider the following relation defined by hearts $\mathcal{H}_{4}\unlhd\mathcal{H}'_{4}$ associated with suitable Gelfand-Tsetlin arrays.

 \begin{centering}

$
  \begin{array}{ccc}
   
     x& y&\\
     &z&        
  \end{array}\unlhd \begin{array}{ccc}
   
     x'& y'&\\
     &z'&        
  \end{array}$ \par\bigskip

  \end{centering}

Then if $x=x', y=y', z'=z-1$ or $x=x', y'=y+1, z'=z$  and the number of Gelfand-Tsetlin arrays is given by a sequence with the form $S_{gt}(n_{1}, r)=n_{1}q_{1}+\underset{j=2}{\overset{2r+1}{\sum}}q_{j}$ then $\mathcal{H}'_{4}$ covers $\mathcal{H}_{4}$ and the associated sum to $\mathcal{H}'_{4}$ has the form $S_{gt}(n'_{1}, r)=n'_{1}q_{1}+\underset{j=2}{\overset{2r+1}{\sum}}q_{j}$ with $n'_{1}=n_{1}-1$. Note that, if $\mathcal{H}_{4}$ and $\mathcal{H}'_{4}$ are incomparable then $x\neq x'$.\hspace{0.5cm}$\square$
\par\bigskip
The next diagram illustrates poset $(\mathscr{H}_{(4,1)},\unlhd),$ where points with the same color are equivalent:

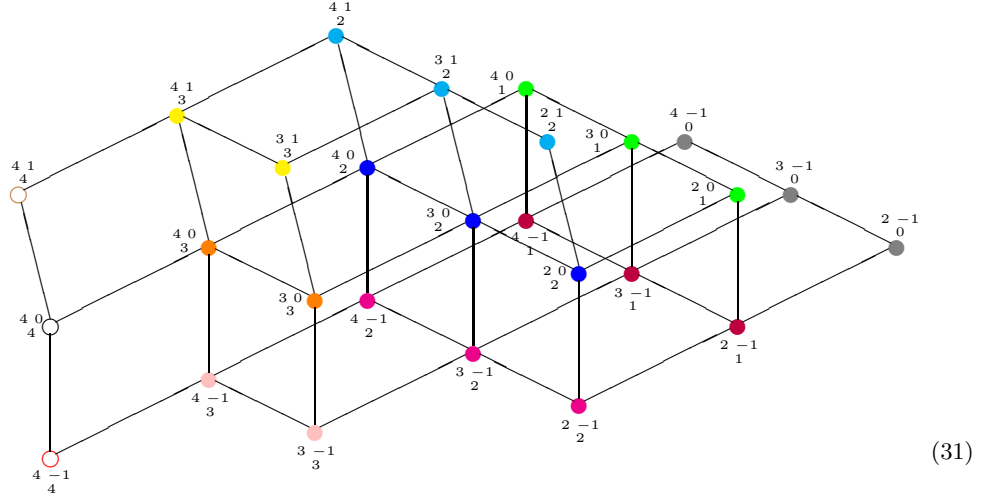
\setlength{\unitlength}{1pt}
\begin{picture}(300,200)
\color{red}\put(0,0){\circle{6}}
\color{pink}\put(60,30){\circle*{6}}
\color{magenta}\put(120,60){\circle*{6}}
\color{purple}\put(180,90){\circle*{6}}
\color{gray}\put(240,120){\circle*{6}}

\color{pink}\put(100,10){\circle*{6}}
\color{magenta}\put(160,40){\circle*{6}}
\color{purple}\put(220,70){\circle*{6}}
\color{gray}\put(280,100){\circle*{6}}

\color{magenta}\put(200,20){\circle*{6}}
\color{purple}\put(260,50){\circle*{6}}
\color{gray}\put(320,80){\circle*{6}}

\color{black}\multiput(2,2)(60,30){4}{\line(2,1){55}}
\color{black}\multiput(102,12)(60,30){3}{\line(2,1){55}}
\color{black}\multiput(202,22)(60,30){2}{\line(2,1){55}}

\color{black}\multiput(98,12)(60,30){4}{\line(-2,1){35.5}}
\color{black}\multiput(198,22)(60,30){3}{\line(-2,1){35.5}}

\color{black}\put(313,87){\tiny{${2\ -1}\atop{0}$}}
\color{black}\put(252,39){\tiny{${2\ -1}\atop{1}$}}
\color{black}\put(192,9){\tiny{${2\ -1}\atop{2}$}}

\color{black}\put(273,107){\tiny{${3\ -1}\atop{0}$}}
\color{black}\put(212,59){\tiny{${3\ -1}\atop{1}$}}
\color{black}\put(152,29){\tiny{${3\ -1}\atop{2}$}}
\color{black}\put(92,-1){\tiny{${3\ -1}\atop{3}$}}

\color{black}\put(233,127){\tiny{${4\ -1}\atop{0}$}}
\color{black}\put(173,80){\tiny{${4\ -1}\atop{1}$}}
\color{black}\put(112,49){\tiny{${4\ -1}\atop{2}$}}
\color{black}\put(52,19){\tiny{${4\ -1}\atop{3}$}}
\color{black}\put(-8,-10){\tiny{${4\ -1}\atop{4}$}}


\color{black}\put(0,50){\circle{6}}
\color{orange}\put(60,80){\circle*{6}}
\color{blue}\put(120,110){\circle*{6}}
\color{green}\put(180,140){\circle*{6}}

\color{orange}\put(100,60){\circle*{6}}
\color{blue}\put(160,90){\circle*{6}}
\color{green}\put(220,120){\circle*{6}}

\color{blue}\put(200,70){\circle*{6}}
\color{green}\put(260,100){\circle*{6}}

\color{black}\multiput(2,52)(60,30){3}{\line(2,1){55}}
\color{black}\multiput(102,62)(60,30){2}{\line(2,1){55}}
\color{black}\multiput(202,72)(60,30){1}{\line(2,1){55}}

\color{black}\multiput(98,62)(60,30){3}{\line(-2,1){35.5}}
\color{black}\multiput(198,72)(60,30){2}{\line(-2,1){35.5}}

\color{black}\multiput(0,3)(60,30){4}{\line(0,1){44.5}}
\color{black}\multiput(100,13)(60,30){3}{\line(0,1){44.5}}
\color{black}\multiput(200,23)(60,30){2}{\line(0,1){44.5}}

\color{black}\put(241,99){\tiny{${2\ 0}\atop{1}$}}
\color{black}\put(185,67){\tiny{${2\ 0}\atop{2}$}}

\color{black}\put(201,119){\tiny{${3\ 0}\atop{1}$}}
\color{black}\put(141,89){\tiny{${3\ 0}\atop{2}$}}
\color{black}\put(85,57){\tiny{${3\ 0}\atop{3}$}}

\color{black}\put(165,141){\tiny{${4\ 0}\atop{1}$}}
\color{black}\put(105,111){\tiny{${4\ 0}\atop{2}$}}
\color{black}\put(45,81){\tiny{${4\ 0}\atop{3}$}}
\color{black}\put(-13,49){\tiny{${4\ 0}\atop{4}$}}

\color{brown}\put(-12,100){\circle{6}}
\color{yellow}\put(48,130){\circle*{6}}
\color{cyan}\put(108,160){\circle*{6}}

\color{yellow}\put(88,110){\circle*{6}}
\color{cyan}\put(148,140){\circle*{6}}

\color{cyan}\put(188,120){\circle*{6}}

\color{black}\multiput(-10,102)(60,30){2}{\line(2,1){55}}
\color{black}\multiput(90,112)(60,30){1}{\line(2,1){55}}

\color{black}\multiput(86,112)(60,30){2}{\line(-2,1){35.5}}
\color{black}\multiput(186,122)(60,30){1}{\line(-2,1){35.5}}

\color{black}\multiput(0,53)(60,30){3}{\line(-1,4){11}}
\color{black}\multiput(100,63)(60,30){2}{\line(-1,4){11}}
\color{black}\multiput(200,73)(60,30){1}{\line(-1,4){11}}

\color{black}\put(184,127){\tiny{${2\ 1}\atop{2}$}}

\color{black}\put(144,147){\tiny{${3\ 1}\atop{2}$}}
\color{black}\put(84,117){\tiny{${3\ 1}\atop{3}$}}

\color{black}\put(104,167){\tiny{${4\ 1}\atop{2}$}}
\color{black}\put(44,137){\tiny{${4\ 1}\atop{3}$}}
\color{black}\put(-16,107){\tiny{${4\ 1}\atop{4}$}}
\color{black}\put(333,0){(31)}

\end{picture}
\par\bigskip
\par\bigskip
\par\bigskip

\addtocounter{corol}{11}
\begin{corol}
For $r\geq 2$ and $n=4$, the number of facets $\mathscr{F}_{(n,r)}$ in the marked order polytope $Q_{\mathscr{H}_{(n, r)}}$ associated with $\mathscr{H}_{(n, r)}$ is $r(r+1)(3r+2).$
\end{corol}

\textbf{Proof.} If $y_{0}$ is the largest integer associated with the hearts $
  \begin{array}{ccc}
   
     x& y_{0}&\\
     &z&        
  \end{array}$  of $\mathscr{H}_{(n, r)}$ then subposet $\mathscr{H}_{(n, r, x_{0})}$, which contains all the hearts with this shape is equal to the poset $\mathscr{P}_{\mathrm{sl}_{r}}$. Thus, the number of facets $\mathscr{F}_{(n, r, x_{0})}$ in the corresponding order polytope equals $r(r+1)$. Note that, $\mathscr{H}_{(n, r)}=\underset{j=1}{\overset{r}{\cup}}\mathscr{H}_{(n, r, y_{0}-j)}$, $\mathscr{H}_{(n, r, s)}\cap\mathscr{H}_{(n, r, t)}=\varnothing$, if $s\neq t$, and $\mathscr{F}_{(n, r, y_{0}-j)}=r(r+1)+t_{j+r}+t_{r}-t_{j-1}$,  where $t_{i}$ denotes the $i$th triangular number. Therefore, $\mathscr{F}_{(n, r)}=\underset{j=0}{\overset{r}{\sum}}\mathscr{F}_{(n, r, y_{0}-j)}+$ the number of covers $\mathscr{C}(y, y+1)$ of the form
$
  \begin{array}{ccc}
   
     x& y&\\
     &z&        
  \end{array}\unlhd \begin{array}{ccc}
   
     x& y+1&\\
     &z&        
  \end{array}$ since $\mathscr{C}(y,y+1)=\frac{r(r+1)(2r+1)}{2}$ then $\mathscr{F}_{(n, r)}=\frac{r(r+1)(4r+3)}{2}+\frac{r(r+1)(2r+1)}{2}$. Therefore, $\mathscr{F}_{(n, r)}=r(r+1)(3r+2)$ as desired.\hspace{0.5cm}$\square$

   \subsubsection{Gelfand-Tsetlin Equations}

Gelfand-Tsetlin equations are diophantine equations of the form:
\addtocounter{equation}{1}
\begin{equation}\label{gt(n)}
\begin{split}
gt(3)&=4x_{1}+12x_{2}+8x_{3}=d_{3},\\
gt(n)&=\underset{i=-1}{\overset{n-3}{\sum}}(s_{(n+i)}-t_{(n-2)})x_{(i+2)}+(t_{(2n-3)}-t_{(n-2)})x_{n}=d_{n}, \quad\text{if}\hspace{0.1cm}n\geq4,
\end{split}
\end{equation}
\par\bigskip
where $d$ is a positive integer and $s_{m}$ ($t_{m}$) denotes the $m$th square (triangular) number. Note that, if $d=64$ then $x_{1}=5$, $x_{2}=3$, and $x_{3}=1$ is a solution of $gt(3)$, such numbers $x_{i}$ gives the number of some standard Gelfand-Tsetlin arrays for which $n=4$ and $r=1$ (see Theorem \ref{NGT}). Numbers $d_{i}$, which correspond to solutions of a given equation $gt(r)$ are said to be \textit{Gelfand-Tsetlin numbers} of type $(4,r)$. Note that, there is not a Frobenius number associated with the equation $gt(3)$.

\par\bigskip

\begin{corol}
\textit{For $r>1$ and $n=4$ sums  $S_{gt}(n_{1},r)$ (see identities $\mathrm{(\ref{sequence})}$ and $\mathrm{(\ref{sequence1})}$) associated with points $\mathcal{H}_{n}\in(\mathscr{H}_{(n, r)},\unlhd)$ are Gelfand-Tsetlin numbers of type $(n, r)$.}
\end{corol}

\textbf{Proof.} If $d_{r}=(r+1)^6$ then Propositions \ref{NGT}, \ref{NGT1}, and \ref{NGT2} allows us to establish identities of type $x_{i}=q_{i}$, for $1\leq i\leq 2n+1$. Thus, $gt(r)=d_{r}$ (see equations (\ref{gt(n)})).\hspace{0.5cm}$\square$
\par\bigskip
The next table shows Frobenius numbers associated with some Gelfand-Tsetlin equations:
\par\bigskip
\begin{centering}
$\begin{array}{|c|c|} \hline
\mathrm{Gelfand-Tsetlin}\hspace{0.1cm}\mathrm{equation} & \mathrm{Frobenius\hspace{0.1cm}number} \\\hline
gt(4) & 33 \\\hline gt(5) & 56 \\\hline gt(6) & 133 \\\hline gt(7) & 179 \\\hline
gt(8)& 181\\\hline gt(9)&299\\\hline gt(10)&394\\\hline gt(11)& 535\\\hline
gt(12)&\infty\\\hline 
\end{array}$ \par\bigskip
\end{centering}

\par\bigskip
We let $\mathscr{H}_{(4 ,r, y_{0})}$ denote the subposet of $(\mathscr{H}_{(4,r)},\unlhd)$ consisting of hearts with a fixed entry $\lambda_{(3,2)}=y_{0}$, with the shape $
  \begin{array}{ccc}
   
     x& &y_{0}\\
     &z&        
  \end{array}$, the corresponding equivalence classes have associated a unique number of the form $S_{gt}(n_{1},j)$. If it is assumed that $n_{1}=1$ for any of these numbers then the following result takes place:
  \addtocounter{teor}{4}
  \begin{teor}
  \textit{Numbers $S_{gt}(n_{1},r)$ with $n_{1}=1$ associated with equivalence classes of points of the subposet $\mathscr{H}_{(4, r, y_{0})}$ define a Brauer configuration algebra $\Lambda_{\Gamma_{gt}}$ with length grading induced from the path algebra $\mathbb{F}Q_{\Gamma_{gt}}$.}
\end{teor}

\par\bigskip

\textbf{Proof.} If $n_{1}=1$ for any number $S_{gt}(n_{1},r)$ then such an assignation defines the Brauer configuration $\Gamma_{gt}=(\Gamma_{gt0},\Gamma_{gt1},\mu,\mathcal{O})$, such that:
\begin{equation}
\begin{split}
\Gamma_{gt0}&=\{q_{1},q_{2},\dots, q_{2r+1}\},\hspace{0.1cm}(\text{see identities}\hspace{0.1cm}(\ref{sequence}), (\ref{sequence1})),\\
\Gamma_{gt1}&=\{w(P_{1})=q^{(1)}_{1}q^{(1)}_{2}\dots q^{(1)}_{2r+1}=w(P_{j}), 2\leq j\leq 2r+1\},\\
\mu(q_{i})&=1,\quad 1\leq i\leq 2r+1,\\
P_{i}&<P_{i+1},\quad 1\leq i\leq 2r\quad (\text{for the construction of successor sequences}),\\
\end{split}
\end{equation}
 where $w(P_{i})$ denotes the word associated with the polygon $P_{i}$ (see (\ref{word})). In this case the Brauer quiver has the form:
 
 \begin{equation}\label{mutfinal}
Q_{\Gamma_{gt}}=
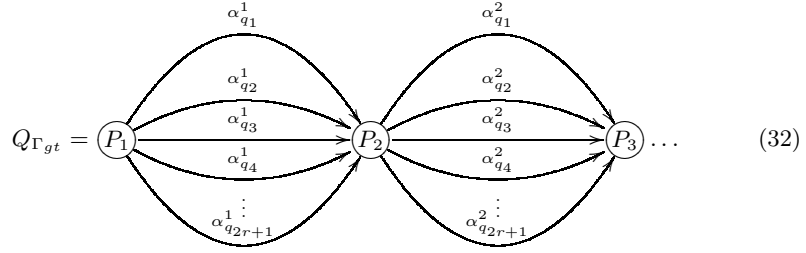
\xymatrix@=80pt{*+[o][F-]{P_{1}}\ar@/^0pt/@[\blue][r]^{\alpha^{1}_{q_{3}}}\ar@/^15pt/@[\blue][r]^{\alpha^{1}_{q_{2}}}\ar@/_15pt/@[\blue][r]^{\alpha^{1}_{q_{4}}}_{\vdots}\ar@/^40pt/@[\blue][r]^{\alpha^1_{q_{1}}}\ar@/_40pt/@[\blue][r]^{\alpha^1_{q_{2r+1}}}&*+[o][F-]{P_{2}}\ar@/^0pt/@[\blue][r]^{\alpha^{2}_{q_{3}}}\ar@/^15pt/@[\blue][r]^{\alpha^{2}_{q_{2}}}\ar@/_15pt/@[\blue][r]^{\alpha^{2}_{q_{4}}}_{\vdots}\ar@/^40pt/@[\blue][r]^{\alpha^2_{q_{1}}}\ar@/_40pt/@[\blue][r]^{\alpha^2_{q_{2r+1}}}&*+[o][F-]{P_{3}}}\dots
\end{equation}
 \par\bigskip
 The admissible ideal is generated via equivalence of special cycles, products of the form $\alpha^{i}_{q_{j}}\alpha^{i+1}_{q_{j'}}$ with $i\neq j'$, $\underset{j=1}{\overset{2r+1}{\prod}}\alpha^{j}_{q_{s}}\alpha^{2r+2}_{q_{s}}\alpha^{1}_{q_{s}}$, where $\underset{j=1}{\overset{2r+1}{\prod}}\alpha^{j}_{q_{s}}\alpha^{2r+2}_{q_{s}}$ is a special cycle and $\alpha^{2r+2}_{q_{s}}$ is an arrow connecting polygons $P_{2r+1}$ and $P_{1}$.  Since $\mu(q_{i})val(q_{i})=2r+1$, for any $q_{i}\in \Gamma_{gt0}$, the result follows as a consequence of Theorem \ref{grading}.\hspace{0.5cm}$\square$

\section{On Diophantine Equations of Type $\mathcal{D}(n_{1},n_{2},\mathcal{K}_{m})$}

In this section, we give criteria to solve equations of type $\mathcal{D}(n_{1},n_{2},\mathcal{K}_{m})$, i.e., diophantine equations of the form:

\begin{equation}\label{D1}
\begin{split}
\underset{i=1}{\overset{m}{\sum}}\lambda_{i}&=n_{1},\\
\underset{j=1}{\overset{m}{\sum}}k_{j}\lambda_{j}&=n_{2}
\end{split}
\end{equation}

with fixed $k_{1},k_{2},\dots, k_{m}$, $n_{1}$ and $n_{2}$.\par\bigskip
Regarding equations of type (\ref{D1}), we have the following results:

\begin{teor}
If $\lambda_{1},\lambda_{2},\dots,\lambda_{m}$ is solution of $\mathcal{D}(n_{1},n_{2},\textbf{m}=\{1,2,\dots, m\})$ and $k_{i}=i$ then $(\lambda_{m},\lambda_{m-1},\dots,\lambda_{1})$ is solution of the diophantine equation $\mathcal{D}(n_{1},(m+1)n_{1}-n_{2},\textbf{m})$.

\end{teor}

\textbf{Proof.} If $(\lambda_{1},\lambda_{2},\dots,\lambda_{m})$ is solution of $\mathcal{D}(n_{1},n_{2},\textbf{m})$ then $\underset{i=1}{\overset{m}{\sum}}\lambda_{i}=n_{1}$ and $\underset{i=1}{\overset{m}{\sum}}i\lambda_{i}=n_{2}$. Thus, 

\begin{equation}
(m+1)n_{1}-n_{2}=(m+1)\underset{i=1}{\overset{m}{\sum}}\lambda_{i}-\underset{i=1}{\overset{m}{\sum}}i\lambda_{i}=\underset{i=1}{\overset{m}{\sum}}(m+1-i)\lambda_{i}=\underset{i=1}{\overset{m}{\sum}}i\lambda_{m+1-i}.
\end{equation}

Therefore $(\lambda_{m},\lambda_{m-1},\dots,\lambda_{1})$ is solution of the desired equation. \par\bigskip Conversely, if $(\lambda_{m},\lambda_{m-1},\dots,\lambda_{1})$ is solution of $\mathcal{D}(n_{1},(m+1)n_{1}-n_{2},\textbf{m})$ then the result holds provided that $(m+1)n_{1}-n'_{2}=n_{2}$, if $n'_{2}=(m+1)n_{1}-n_{2}$.\hspace{0.5cm}$\square$

\begin{teor}
An equation of type $\mathcal{D}(n_{1},n_{2},\textbf{m})$ with $k_{i}=i$ has at least one solution if and only if $n_{1}+\frac{m(m-1)}{2}\leq n_{2}\leq mn_{1}-\frac{m(m+1)}{2}$.
\end{teor}

\textbf{Proof.} If $x=(\lambda_{1},\lambda_{2},\dots,\lambda_{m})\in\mathbb{N}^{m}$ and $\lambda_{i}<\lambda_{j}$ for $i<j$ then for
\begin{equation}
x'=(\lambda_{1},\dots,\lambda_{j},\dots,\lambda_{i},\dots,\lambda_{m})\notag
\end{equation}
 it holds that $\underset{i=1}{\overset{m}{\sum}}ix_{i}<\underset{i=1}{\overset{m}{\sum}}ix'_{i}$, where $x_{i}$ ($x'_{i}$) denotes the $i$th coordinate of the vector $x_{i}$ ($x'_{i}$), in general $\underset{i=1}{\overset{m}{\sum}}x_{i}\lambda_{i}\leq \underset{i=1}{\overset{m}{\sum}}x'_{i}\lambda_{i}$, if $\lambda_{1}\leq \lambda_{2}\leq\dots\leq\lambda_{m}$. Therefore, the minimum value of $n_{2}$ for which equation $\mathcal{D}(n_{1},n_{2},\textbf{m})$ has a solution is attained at $x=(n_{1}-m+1, 1,1,\dots, 1)$ whereas the maximum  value of $n_{2}$ satisfying the condition is attained at $x'=(1,1,\dots,1,n_{1}-m+1)$. Thus equation $\mathcal{D}(n_{1},n_{2},\textbf{m})$ has a solution if the following conditions for $n_{1}$ and $n_{2}$ hold:
\begin{equation}
\begin{split}
\underset{i=1}{\overset{m}{\sum}}ix_{i}&\leq n_{2}\leq \underset{i=1}{\overset{m}{\sum}}ix'_{i},\\
n_{1}-m+\underset{i=1}{\overset{m}{\sum}}i&\leq n_{2}\leq \underset{i=1}{\overset{m-1}{\sum}}i+m(n_{1}-m).\notag
\end{split}
\end{equation}

Therefore, 
\begin{equation}
n_{1}-m+\frac{m(m+1)}{2}\leq n_{2}\leq \underset{i=1}{\overset{m-1}{\sum}}i+m(n_{1}-m).\notag
\end{equation} 
 
 And
 
 \begin{equation}
 n_{1}+\frac{m(m-1)}{2}\leq n_{2}\leq mn_{1}-\frac{m(m+1)}{2}.\hspace{0.5cm}\notag\square
 \end{equation}

\begin{teor}\label{main}
For a fixed integer positive $m_{0}$, $n\geq2$, and a fixed seed $(\Gamma,\mathscr{X})$ (see $\mathrm{(\ref{mut01})}$ and $\mathrm{(\ref{mut02})}$), the Brauer configuration algebra $\Lambda_{\Phi^{m_{0}}}$ induced by the Brauer configuration 
\begin{equation}
\Phi^{m_{0}}=(\Phi^{m_{0}}_{0},\Phi^{m_{0}}_{1},\mu,\mathcal{O})\quad (see\hspace{0.1cm} \mathrm(\ref{mutation}))\notag
\end{equation}

has associated a finite non-deterministic automaton whose states are given by solutions of a problem of type $\mathcal{D}(n_{1},l^{2}2^{n-2},\mathcal{K}_{m_{i}})$ with $n_{1}\leq 16$, $l=4t$, for some $t\geq1$, and 
\begin{equation}
\mathcal{K}_{m_{i}}= \{\mathrm{occ}(T(\alpha), T(\mathcal{M}^{i}))\mid T(\alpha)\in T(\Phi^{m_{0}}_{0}), T(\mathcal{M}^{i})\in T(\Phi^{m_{0}}_{1}), \alpha\in\Phi^{m_{0}}_{0}, \mathcal{M}^{i}\in \Phi^{m_{0}}_{1} \}\notag
\end{equation}

where $T(\Phi^{m_{0}})=((T(\Phi^{m_{0}}_{0}),T(\Phi^{m_{0}}_{1}),T(\mu),T(\mathcal{O}))$ is a suitable transformation between Brauer configurations.
\end{teor}

\textbf{Proof.} Since the length of the message $|M(\Phi^{m_{0}})|=l^{2}2^{n}$ then it consists of $l^{2}2^{n-2}$ lists of four bits. We let $\mathscr{A}=\{\alpha_{1},\alpha_{2},\dots,\alpha_{l^{2}2^{n-2}}\}$ denote this set of lists. Define a map $T:\mathscr{A}\rightarrow Hex$ such that $T(\alpha_{i})=\alpha'_{i}\in Hex$, where $Hex$ is a notation for the hexadecimal numbering system $Hex=\{0,1,2,\dots,9,A,B,C,\dots, F\}$. In this case, $T(\alpha_{i}\alpha_{i+1})=\alpha'_{i}\alpha'_{i+1}\in Hex^{*}$ if $\alpha_{i},\alpha_{i+1}\in\mathscr{A}$.\par\bigskip
The map $T$ defines a new Brauer configuration $T(\Phi^{m_{0}})=(\mathrm{Img}\hspace{0.1cm}T, T(\Phi^{m_{0}}_{1}),\mu',\mathcal{O}')$ which we assume reduced without loss of generality. Each polygon $T(\mathcal{M}^{i})$ consists of elements of the form $\alpha'_{j}$ with $\alpha_{j}\in\mathcal{M}^{i}$. Actually, if $M(\mathcal{M}^{i})=\alpha_{1}\alpha_{2}\dots \alpha_{r_{i}}$ is the message of $\mathcal{M}^{i}$ then $\alpha'_{1}\alpha'_{2}\dots\alpha'_{r_{i}}\in \mathscr{A}^{*}$ is the message of $T(\mathcal{M}^{i})$. Besides, if $\mathcal{M}_{i}<\mathcal{M}_{i+1}$ in $\Phi^{m_{0}}$ then $T(\mathcal{M}_{i})<T(\mathcal{M}_{i+1})$ in $T(\Phi^{m_{0}})$. And $\mu'(\alpha'_{i})=1$ for any $\alpha'_{i}\in \mathrm{Img}\hspace{0.1cm}T$. Thus, the message $M(T(\Phi^{m_{0}}))\in\mathscr{A}^{*}$ is a $l^{2}2^{n-2}$ word whose letters $\alpha'_{j}$ can be grouped according to its valency. Thus, $M(T(\Phi^{m_{0}}))$ has the form:

\begin{equation}
M(T(\Phi^{m_{0}}))=\mathscr{A}_{i_{1}}\mathscr{A}_{i_{2}}\dots\mathscr{A}_{i_{m}},
\end{equation}

where $\mathscr{A}_{i_{s}}$ is a multiset with $|\mathscr{A}_{i_{s}}|=L_{i_{s}}$ and $\mathscr{A}_{i_{x}}\cap\mathscr{A}_{i_{y}}=\varnothing$. Note that $\mathscr{A}_{i_{s}}$ consists of all letters $\alpha'_{i}$ such that $val(\alpha'_{i})=v_{i_{s}}$, i.e., the message $M(\mathscr{A}_{i_{s}})$ associated with $\mathscr{A}_{i_{s}}$ can be written as 
\begin{equation}
M(\mathscr{A}_{i_{s}})=(\alpha'_{is_{1}}\alpha'_{is_{2}}\dots\alpha'_{L_{i_{s}}})^{v_{i_{s}}}.
\end{equation}

Therefore, 

\begin{equation}
\begin{split}
\underset{h=1}{\overset{m}{\sum}}L_{i_{h}}&=n_{1}\leq16,\\
\underset{g=1}{\overset{m}{\sum}}L_{i_{g}}v_{i_{g}}&=l^{2}2^{n-2}.
\end{split}
\end{equation}

Then terms $L_{i_{h}}$ give a solution of a diophantine equation of type
\begin{equation}
 \mathcal{D}(n_{1},l^{2}2^{n-2},\mathcal{K}_{m}=\{v_{i_{1}},\dots, v_{i_{m}}\}).
 \end{equation}
 
Since any Brauer configuration algebra defines a regular language whose associated automaton uses arrows of the Brauer quiver as transitions between states given by polygons, which in this case are obtained by mutation, thus the probability $P_{r}(M)$ that a given message $M$ or diophantine equation occurs after applying a mutation to a fixed seed is such that $0<P_{r}(M)<1$ (see Remark \ref{probability}), we conclude that the associated automaton is non deterministic.\hspace{0.5cm}$\square$
 \par\bigskip
 As an example suppose that a seed $(\Gamma,\mathscr{X})$ is defined by the polynomial $p(x)=x^{8}+x^{4}+x^{3}+1$, and that we use set $Hex=\{0=0000,1=0001,2=0010,\dots,A=1010,\dots,F=1111\}$ to denote polynomials. Then $\mathscr{X}=(x_{1},x_{2},x_{3},x_{4})$ can be denoted as follows:
\begin{equation}
\begin{split}
x_{1}&=(AF,C0,13,10),\\
x_{2}&=(D0,B3,8A,F2),\\\notag
x_{3}&=(CE,C4,62,3D),\\
x_{4}&=(A2,74,79,7D).
\end{split}
\end{equation}

 \begin{equation}\label{subbytes0}
 \begin{split}
 \Gamma&=(\Gamma_{0},\Gamma_{1},\mu,\mathcal{O}),\\
 \Gamma_{0}&=\{0,1\},\\
 \Gamma_{1}&= \{w(0)=\{10101111\dots00010000\},w(1)=\{11010000\dots11110010\},\dots\},\\
 \mu(0)&=\mu(1)=1,\\
 w(0)&<w(1)<w(2)<w(3).\\
  \end{split}
 \end{equation}
 
 Then 
 
 \begin{equation}
M(\Gamma)=(AC03D27)^{(3)}(F14)^{(2)}(B8E69)^{(1)} \notag
 \end{equation}
 \par\bigskip
 which builds a solution for a diophantine equation $\mathcal{D}(15,32,\{3,2,1\})$ with the form:
 
 \begin{equation}
 \begin{split}
 \lambda_{1}+\lambda_{2}+\lambda_{3}&=15,\\
 3\lambda_{1}+2\lambda_{2}+\lambda_{3}&=32.\notag
 \end{split}
 \end{equation}
 
 For a given polynomial $p(x)=w(i)=a_{0}+a_{1}x+a_{2}x^{2}+a_{3}x^{3}+a_{4}x^{4}+a_{5}x^{5}+a_{6}x^{6}+a_{7}x^{7}\in\mathbb{F}$  the map $\tau(p(x))=\tau(x_{i})$ associated with a mutation of the seed $(\Gamma,\mathscr{X})$ is defined in such a way that
 
\begin{equation}\label{subbytes1}
\begin{split}
\tau(x_{i})&=\underset{s=0}{\overset{7}{\sum}}b_{s}x^{s}+v_{i/4,0},\\
b_{j}&=a_{j}+a_{j+4}+a_{j+5}+a_{j+6}+a_{j+7}+c_{j},\\
(c_{0},c_{1},c_{2},c_{3},c_{4},c_{5},c_{6},c_{7})&=(1,1,0,0,0,1,1,0),\\
\mathscr{H}(x_{1},x_{2},x_{3},x_{4})&=(\tau(x_{2}),\tau(x_{3}),\tau(x_{4}),\tau(x_{1})),\quad m_{0}=10.
\end{split}
\end{equation} 
 
The following is the list of vectors $v_{j,0}$, $1\leq j\leq 10$. 

\begin{equation}\label{subbytes2}
\begin{split}
v_{1,0}&\leftarrow01000000\\
v_{2,0}&\leftarrow02000000\\
v_{3,0}&\leftarrow04000000\\
v_{4,0}&\leftarrow08000000\\
v_{5,0}&\leftarrow10000000\\
v_{6,0}&\leftarrow20000000\\
v_{7,0}&\leftarrow41000000\\
v_{8,0}&\leftarrow81000000\\
v_{9,0}&\leftarrow1B000000\\
v_{10,0}&\leftarrow36000000
\end{split}
\end{equation}
 
 For $m_{0}=10$, any mutation of this seed gives rise to a solution of a diophantine equation of type $(n_{1}\leq 16,32,\mathcal{K}_{m})$ with $\mathcal{K}_{m}$ being a set of the form:
 
 \begin{equation}
\begin{split} 
\mathcal{K}_{m}&=\{1,2,3\},\\
\mathcal{K}_{m}&=\{1,2,3,4\},\\
\mathcal{K}_{m}&=\{1,2,3,4,5\},\\\notag
\mathcal{K}_{m}&=\{1,2,3,4,7\},\\
\mathcal{K}_{m}&=\{1,2,3,4,7\},\\
\mathcal{K}_{m}&=\{1,2,3,4,8\}. 
\end{split} 
 \end{equation}
 \par\bigskip
\begin{centering}
 \textbf{A Cryptographic Application}\par\bigskip
 \end{centering}
 
 As an application of the notion of mutation of a Brauer configuration, in this section, we present the schedule of an AES key (Advanced Encryption Standard) or Rijndael which is the cornerstone of wireless communications \cite{Stinson}.\par\bigskip

 AES is a symmetric block cipher also considered as a substitution-permutation network, which was adopted in 2001 by the US government as the current standard cryptographic protocol. It is considered secure against different types of attacks. It requires keys of 128 bits (for ten rounds), 192 bits (for 12 rounds) and 256 bits (for 14 rounds). Encryption and decryption algorithms are carried out via polynomials defined over suitable Galois fields. Many routers provide protocols WPA2-PSK (TKIP), WPA2-PSK (AES), and WPA2-PSK (TKIP/AES) as options to ensure WI-FI security. The new protocol WPA3 uses keys of 128 bits (CCMP-128), 192 bits (WPA3-Enterprise mode) to secure wireless computer networks. \par\bigskip

 In the cryptosystem AES, a plaintext also called state is a sequence of 16 bytes, the encryption process also generates a 16-bytes sequence by using keys of 128, 192, or 256 bits. Such length depends on the number of rounds developed in the encryption process, 10, 12 or 14 respectively.\par\bigskip
 
 Each round in an encryption process requires 4 different transformations:
 
 \begin{enumerate}
\item SubBytes,
\item ShiftRows,
\item MixColumns,
\item AddRoundKey.
 \end{enumerate}
 
 For the last round the function MixColums is not executed.\par\bigskip
 
 The next table gives all the possible outputs of the transformation SubBytes \cite{Stinson}:
 \par\bigskip
 \begin{centering}
 $\scalebox{0.73}{
\begin{tabular}{||c||c|c|c|c|c|c|c|c|c|c|c|c|c|c|c|c|}\hline
&0&1&2&3&4&5&6&7&8&9&A&B&C&D&E&F\\\hline\hline
0&63&7C&77&7B&F2&6B&6F&C5&30&01&67&2B&FE&D7&AB&76\\\hline
1&CA&82&C9&7D&FA&59&47&F0&AD&D4&A2&AF&9C&A4&72&C0\\\hline
2&B7&FD&93&26&36&3F&F7&CC&34&A5&E5&F1&71&D8&31&15\\\hline
3&04&C7&23&C3&18&96&05&9A&07&12&80&E2&EB&27&B2&75\\\hline
4&09&83&2C&1A&1B&6E&5A&A0&52&3B&D6&B3&29&E3&2F&84\\\hline
5&53&D1&00&ED&20&FC&B1&5B&6A&CB&BE&39&4A&4C&58&CF\\\hline
6&D0&EF&AA&FB&43&4D&33&85&45&F9&02&7F&50&3C&9F&A8\\\hline
7&51&A3&40&8F&92&9D&38&F5&BC&B6&DA&21&10&FF&F3&D2\\\hline
8&CD&0C&13&EC&5F&97&44&17&C4&A7&7E&3D&64&5D&19&73\\\hline
9&60&81&4F&DC&22&2A&90&88&46&EE&B8&14&DE&5E&0B&DB\\\hline
A&E0&32&3A&0A&49&06&24&5C&C2&D3&AC&62&91&95&E4&79\\\hline
B&E7&C8&37&6D&8D&D5&4E&A9&6C&56&F4&EA&65&7A&AE&08\\\hline
C&BA&78&25&2E&1C&A6&B4&C6&E8&DD&74&1F&4B&BD&8B&8A\\\hline
D&70&3E&B5&66&48&03&F6&0E&61&35&57&B9&86&C1&1D&9E\\\hline
E&E1&F8&98&11&69&D9&8E&94&9B&1E&87&E9&CE&55&28&DF\\\hline
F&8C&A1&89&0D&BF&E6&42&68&41&99&2D&0F&B0&54&BB&16\\\hline
\end{tabular}

}$

\par\bigskip
\end{centering}

 The key schedule is the process for which all the keys to be used in the encryption process are generated, such keys are called subkeys. For keys of 128 bits (or 16 bytes) of length the process generates 11 subkeys, the initial key, the 9 main rounds and the final round.\par\bigskip
 
 The expanded key can be seen as an array of 32-bit words numbered from 0 to 43 (0 for the initial key), words that are a multiple 4 $(w_{4}, w_{8},\dots, w_{40})$ are calculated as follows:
 
 \begin{enumerate}
 \item Applying the RotWord and SuBytes transformation to the previous word $w_{i-1}$, (SubBytes (RotWord$(x_{1},x_{2},x_{3},x_{4}))=\mathrm{SuBytes}(x_{2},x_{3},x_{4},x_{1}))$.
 \item Adding (XOR) this result to the word 4 positions earlier $w_{i-4}$ plus a round constant called RCON.
 \item The remaining 32-bit words $w_{i}$ are calculated by adding (XOR) the previous word $w_{i-1}$, with the word 4 positions earlier.
 \end{enumerate}
 We let $K^{n,a}$ denote the $n$th subkey of an AES schedule, where the original key $K=K^{0,a}$ has 128 bits. Besides, $S_{n}(K^{0,a})=K^{n+1,a}$ is the result of the AES schedule after the $n$th round.
 \addtocounter{corol}{4}
\begin{corol}
There exists an integer $N$ for which $S_{N-1}(K^{0,a})=K^{M,a}$, for some $0\leq M<N$.

\end{corol}

\textbf{Proof.} For $l=4$, $p(x)=x^{8}+x^{4}+x^{3}+1$ and vectors $\mathrm{RCON}$ given in $\mathrm(\ref{subbytes2})$, it holds that the set of polygons $\Phi^{m_{0}}_{1}$ (see $\mathrm(\ref{mutation})$) of $\Phi^{m_{0}}$ obtained by the mutation rules given in $\mathrm{(\ref{subbytes1})}$ is the $\mathrm{AES}$-key schedule of a seed-key $(\Gamma,\mathscr{X}=(x_{1},x_{2},x_{3},x_{4}))=K^{0,a}$. Thus, the result holds as a direct consequence of Theorem \ref{mutteor}.\hspace{0.5cm}$\square$

\par\bigskip
Agustín Moreno Ca\~{n}adas\\
Department of Mathematics\\
Universidad Nacional de Colombia\\
amorenoca@unal.edu.co\par\bigskip

Juan David Camacho\\
Department of Mathematics\\
Universidad Nacional de Colombia\\
judcamachov@unal.edu.co\par\bigskip

Isa\'{i}as David Mar\'{i}n Gaviria\\
Deparment of Mathematics\\
Universidad Nacional de Colombia\\
imaringa@unal.edu.co


\begin{bibdiv}
\begin{biblist}


\bib{MAO}{article}{title={Brauer configuration algebras for multimedia based encryption and security applications},

Author={M.A.O. Angarita}, Author={A.M. Ca\~{n}adas}, Author={},


Journal={Multimedia Tools and Applications}, number={},  volume={},
date={2021}, pages={1-29}, Note={doi.org/10.1007/s11042-020-10239-3} }




\bib{Fernandez1}{article}{title={Homological ideals as integer specializations of some Brauer
Configuration algebras}, subtitle={}, Author={A.M. Ca\~{n}adas}, Author={P.F.F. Espinosa},
			Author={}, Author={}, Author={}, journal={Ukrainian Journal of Mathematics (accepted)},
			volume={}, date={2020}, number={}, pages ={}, Note={ arXiv:2104.00050v1}}




\bib{Chapman}{article}{title={Factoring in the Chicken McNugget monoid},Author={S.T. Chapman}, Author={C. O'Neill}, Author={},
Journal={Mathematics Magazine}, number={5},  volume={91},
date={2015}, pages={323-336} }





\bib{DeLoera}{article}{title={The many aspects of counting lattice points in polytopes}, subtitle={}, Author={J. De Loera},
			Author={}, Author={}, Author={}, journal={Mathematische semesterberichte},
			volume={52}, date={2005}, number={}, pages ={175-195}}


\bib{Fernandez}{book}{title={Categorification of Integer Sequences and Its Applications}, Author={P.F.F. Espinosa}, Author={}, volume={}, edition={}, date={2020}, Series={}, Publisher={National University of Colombia }, Note={PhD Dissertation} }



\bib{Fisher1}{article}{title={Refined enumerations of alternating sign matrices: $(d,m)$-trapezoids with prescribed top and bottom row}, subtitle={}, Author={I. Fisher},
			Author={}, Author={}, Author={}, journal={J Algebr Combin},
			volume={33}, date={2011}, number={}, pages ={239-257}}






\bib{Fisher2}{article}{title={Sequences of labeled trees related to Gelfand-Tsetlin patterns}, subtitle={}, Author={I. Fisher},
			Author={}, Author={}, Author={}, journal={Advances in Applied Mathematics},
			volume={49}, date={2012}, number={}, pages ={165-195}}








\bib{Fourier}{article}{title={Marked poset polytopes: Minkowski sums, indecomposables, and unimodular equivalence}, subtitle={}, Author={G. Fourier},
			Author={}, Author={}, Author={}, journal={Journal of Pure and Applied Algebra},
			volume={220}, date={2016}, number={2}, pages ={606-620}}
			
			
	\bib{Futorny}{article}{title={Irreducible generic Gelfand-Tsetlin modules}, subtitle={}, Author={V. Futorny},
			Author={D. Grantcharov}, Author={L.E. Ramirez}, Author={}, journal={SIGMA},
			volume={11}, date={2015}, number={}, pages ={239-257}}
			
			

			
	\bib{Garcia}{book}{title={Numerical semigroups},
subtitle={}, Author={P. Garc\'{\i}a-S\'anchez}, Author={J.C. Rosales}, journal={},
publisher={Springer}, volume={}, date={2009}, pages={},
address={New York}, Note={ISBN 978-1-4419-0160-6}

}		
			

\bib{Green}{article}{title={Brauer configuration algebras: A generalization of
				Brauer graph algebras}, subtitle={}, Author={E.L. Green},
			Author={}, Author={S. Schroll}, Author={}, journal={Bull. Sci. Math.},
			volume={141}, date={2017}, number={}, pages ={539--572}}
		











\bib{Ramirez}{article}{title={Complexity of the Frobenius problem},
subtitle={}, Author={J.L. Ram\'{\i}rez-Alfons\'{\i}n}, Author={}, journal={Combinatorica},
volume={16}, date={1996}, pages={143--147}
}


\bib{RamirezF}{article}{title={Combinatorics of Irreducible Gelfand-Tsetlin $sl(3)$-modules},
subtitle={}, Author={L.E. Ramirez}, Author={}, journal={Algebra and Discrete Mathematics},
volume={14}, date={2012}, pages={276--296}
}







\bib{Rees}{article}{title={The Automata that define Representations
of monomial algebras}, Author={S. Rees}, Author={},Author={}, volume={11}, edition={}, date={2008}, Series={},
Journal={Algebr Represent Theor}, Pages={207-214} }



\bib{Ballester}{article}{title={Varieties and covarieties of languages}, Author={J. Rutten}, Author={A. Ballester-Bolinches},Author={E.C. i Ll\'opez}, volume={298}, edition={}, date={2013}, Series={},
Journal={ENTCS}, Pages={7-28} }







\bib{Sierra}{article}{title={The dimension of the center of a Brauer configuration algebra}, subtitle={}, Author={A. Sierra},
			Author={}, Author={}, Author={}, journal={J. Algebra},
			volume={510}, date={2018}, number={}, pages ={289-318}
			
		}
		


\bib{Stinson}{book}{title={Cryptography: Theory and Practice},
subtitle={}, Author={D.R. Stinson}, Author={M. Paterson}, journal={},
publisher={Chapman and Hall/CRC Press}, volume={}, date={2019}, pages={1--598},
address={}, Note={ISBN 9781138197015}

}		

\bib{Zeilberger}{article}{title={Proof of the refined alternating sign matrix conjecture}, subtitle={}, Author={D. Zeilberger},
Author={}, Author={}, Author={}, journal={Electron. J. Combin},
volume={3}, date={1996}, number={2}, pages ={165-195}, Note={paper 13}}
		
		






\end{biblist}
\end{bibdiv}
\end{document}